\documentclass[a4paper,11pt]{amsart}
\usepackage{amssymb,amsmath}
\usepackage{amsopn}
\usepackage{mathrsfs}
\usepackage{color}
\usepackage{mathtools}
\usepackage{wasysym}
\usepackage{vmargin}
\usepackage{pdfsync}
\usepackage{vmargin}
\usepackage{color}
\usepackage{amscd}
\usepackage{verbatim}
\usepackage{bbm}
\usepackage{tikz}
\usetikzlibrary{patterns}

\newtheorem{theorem}{Theorem}[section]
\newtheorem{definition}[theorem]{Definition}
\newtheorem{lemma}[theorem]{Lemma}
\newtheorem{proposition}[theorem]{Proposition}
\newtheorem{exemple}[theorem]{Example}
\newtheorem{corollary}[theorem]{Corollary}
\newtheorem{remark}[theorem]{Remark}
\numberwithin{equation}{section}

\newcommand{\rr}{\mathbb{R}}
\newcommand{\eps}{\varepsilon}
\newcommand{\nn}{\mathbb{N}}

\def\un{{\mathrm{1~\hspace{-1.4ex}l}}}

\def\R{\mathbb R}

\def\l2{L^2(\R^{n})}
\def\L2{L^2(\R^{2n})}
\def\supp{\operatorname{supp}}

\def\eps{\varepsilon}

\def\mat22#1#2#3#4{\begin{pmatrix}#1&#2\\ #3&#4\end{pmatrix}}

\begin{document}
\title[Exact controllability of the free and harmonic Schr\"odinger equation]{Geometric conditions for the exact controllability of fractional free and harmonic Schr\"odinger equations}

\author{J\'er\'emy \textsc{Martin}, Karel \textsc{Pravda-Starov}}

\address{\noindent \textsc{J\'er\'emy Martin, Univ Rennes, CNRS, IRMAR - UMR 6625, F-35000 Rennes, France
}}
\email{jeremy.martin@univ-rennes1.fr}

\address{\noindent \textsc{Karel Pravda-Starov, Univ Rennes, CNRS, IRMAR - UMR 6625, F-35000 Rennes, France
}}
\email{karel.pravda-starov@univ-rennes1.fr}


\begin{abstract}
We provide necessary and sufficient geometric conditions for the exact controllability of the one-dimensional fractional free and fractional harmonic Schr\"odinger equations. The necessary and sufficient condition for the exact controllability of fractional free Schr\"odinger equations is derived from the Logvinenko-Sereda theorem and its quantitative version established by Kovrijkine,  whereas the one for the exact controllability of fractional harmonic Schr\"odinger equations is deduced from an infinite dimensional version of the Hautus test for Hermite functions and the Plancherel-Rotach formula.

\end{abstract}

\maketitle

\section{Introduction}

The purpose of this article is to provide necessary and sufficient geometric conditions on control subsets to ensure exact controllability of fractional free Schr\"odinger equations and fractional harmonic Schr\"odinger equations.
The notion of exact controllability is defined as follows:

\medskip

\begin{definition} [Exact controllability] Let $B$ be a closed operator on $L^2(\rr^d)$ which is the infinitesimal generator of a strongly continuous semigroup $(e^{-tB})_{t \geq 0}$ on $L^2(\rr^d)$, $T>0$ and $\omega$ be a measurable subset of $\mathbb{R}^d$. 
The equation 
\begin{equation}\label{syst_general}
\left\lbrace \begin{array}{ll}
(\partial_t + B)f(t,x)=u(t,x)\un_{\omega}(x)\,, \quad &  x \in \mathbb{R}^d,\ t>0, \\
f|_{t=0}=f_0 \in L^2(\rr^d),                                       &  
\end{array}\right.
\end{equation}
is said to be exactly controllable from the control subset $\omega$ in time $T>0$ if, for any initial datum $f_0 \in L^{2}(\mathbb{R}^d)$ and any final datum $f_T \in L^2(\rr^d)$, there exists $u \in L^2((0,T)\times\mathbb{R}^d)$, supported in $(0,T)\times\omega$, such that the mild (or semigroup) solution of \eqref{syst_general} satisfies $f(T,\cdot)=f_T$. The equation \eqref{syst_general} is said to be exactly controllable from the control subset~$\omega$ if there exists $T>0$ such that this system is exactly controllable from $\omega$ in time $T>0$.
\end{definition}

\medskip

By the Hilbert uniqueness method, see \cite{coron_book, JLL_book}, the exact controllability of the evolution equation \eqref{syst_general} is equivalent to the exact observability of the adjoint system 
\begin{equation} \label{adj_general}
\left\lbrace \begin{array}{ll}
(\partial_t + B^*)g(t,x)=0\,, \quad & x \in \mathbb{R}^d\,,\ t>0 \,, \\
g|_{t=0}=g_0 \in L^2(\rr^d),
\end{array}\right.
\end{equation}
where $B^*$ denotes the $L^2(\rr^d)$-adjoint of the operator $B$.
The notion of exact observability is defined as follows:

\medskip

\begin{definition} [Exact observability] Let $T>0$ and $\omega$ be a measurable subset of $\mathbb{R}^d$. 
The evolution equation \eqref{adj_general} is said to be exactly observable from the control subset $\omega$ in time $T>0$ if there exists a positive constant $C_T>0$ such that,
for any initial datum $g_0 \in L^{2}(\mathbb{R}^d)$, the mild (or semigroup) solution of \eqref{adj_general} satisfies
\begin{equation*}\label{eq:observability}
\int\limits_{\mathbb{R}^d} |g_0(x)|^{2} dx  \leq C_T \int\limits_{0}^{T} \Big(\int\limits_{\omega} |g(t,x)|^{2} dx\Big) dt\,.
\end{equation*}
\end{definition}
In this paper, we study evolution equations associated to skew-selfadjoint operators. More precisely, we consider evolution equations
\begin{equation}\label{obs01}
\left\lbrace \begin{array}{ll}
\partial_tf(t,x) + i(Af)(t,x)=0 \,, \quad &  x \in \mathbb{R}^d,\ t>0, \\
f|_{t=0}=f_0 \in L^2(\rr^d),                                       &  
\end{array}\right.
\end{equation}
where the selfadjoint operators $A$ on $L^2(\rr^d)$ generating a strongly continuous group $(e^{itA})_{t \in \rr}$ on $L^2(\rr^d)$ is specifically given by a fractional Laplacian $A=(-\Delta_x)^s$ equipped with the domain
$$\mathcal{D}\big((-\Delta_x)^s\big)=\big\{f \in L^2(\rr^d) :\ |\xi|^{2s}\widehat{f} \in L^2(\rr_{\xi}^d)\big\},$$
or a fractional harmonic oscillator $A=(- \Delta_x + |x|^2)^s$, with $s>0$, equipped with the domain
$$\mathcal{D}\big((- \Delta_x + |x|^2)^s\big)=\Big\{f \in L^2(\rr^d) :\ \sum_{\alpha \in \nn^d}(2|\alpha|+d)^{2s}|\langle f,\Psi_{\alpha}\rangle_{L^2(\rr^d)}|^2 <+\infty\Big\},$$
where $(\Psi_{\alpha})_{\alpha \in \nn^d}$ denotes the Hermite basis of $L^2(\rr^d)$ defined in (\ref{e1}) and (\ref{bbb2}).

Our starting point is the following result proved by Miller. This author shows in \cite[Corollary 2.17]{miller} that the observability of the system (\ref{obs01}) is equivalent to the following spectral estimates:

\medskip

\begin{proposition}[Miller {\cite[Corollary~2.17]{miller}}]\label{ineqspec}
Let $A$ be a selfadjoint operator on $L^2(\rr^d)$, which is the infinitesimal generator of a strongly continuous group $(e^{itA})_{t \in \rr}$ on $L^2(\rr^d)$. If the evolution equation \emph{(\ref{obs01})} is exactly observable from a measurable subset $\omega \subset \rr^d$ at some time $T>0$ then there exist some positive constants $k>0$ and $ D>0$ such that 
\begin{equation}\label{Spectral_ineq}
\forall \lambda \in \rr, \forall f \in \un_{\{|A-\lambda| \leq \sqrt{D}\}} \big(L^2(\rr^d)\big), \quad  \|f\|_{L^2(\rr^d)} \leq \sqrt{k} \|f\|_{L^2(\omega)}.
\end{equation}
Conversely, when the spectral estimates \emph{(\ref{Spectral_ineq})} hold for some $k>0$ and $D>0$, then the system \emph{(\ref{obs01})} is exactly observable from $\omega$ at any time 
\begin{equation}\label{bbb1}
T>\pi \sqrt{\frac{1+k}{D}}.
\end{equation}
\end{proposition}

\medskip

The operators $\un_{\{|A-\lambda| \leq \sqrt{D}\}}$ appearing in \eqref{Spectral_ineq} are defined by functional calculus. It is important to point out that the constants $D$ and $k$ do not depend on $\lambda$.
We recall the proof given by Miller in Appendix (Section~\ref{proof_ineqspec}) since the lower bound obtained by Miller in \cite[Corollary 2.17]{miller} is different to ours in (\ref{bbb1}).

When the operator $A= (- \Delta_x)^s$, with $s >0$, is a fractional Laplacian, the family $(\un_{\{|A-\lambda| \leq \sqrt{D}\}})_{\lambda \in \rr}$ is a family of frequency cutoff operators on bounded frequency subsets and for all $D>0$ and $\lambda \in \rr$,
\begin{equation}\label{cutoff_operator}
 \un_{\{|A-\lambda| \leq \sqrt{D}\}} \big(L^2(\rr^d)\big) = \Big\{ f \in L^2(\rr^d) : \ \supp \hat{f} \subset \big\{ \xi \in \rr^d, \, ||\xi|^{2s}- \lambda | \leq \sqrt{D} \big\} \Big\},
\end{equation}
where $|\cdot|$ denotes the Euclidean norm on $\rr^d$.
The spectral estimates \eqref{Spectral_ineq} are in this case related to the notion of annihilating pairs: 
\medskip

\begin{definition} [Annihilating pairs]
Let $S,\Sigma$ be two measurable subsets of $\rr^d$. 
\begin{itemize}
\item[-] The pair $(S,\Sigma)$ is said to be a weak annihilating pair if the only function $f\in L^2(\rr^d)$
with $\supp f\subset S$ and $\supp\widehat{f} \subset \Sigma$ is zero $f=0$.
\item[-] The pair $(S,\Sigma)$ is said to be a strong annihilating pair if there exists a positive constant $C=C(S,\Sigma)>0$ such that 
for all $f \in L^2(\rr^d)$,
\begin{equation*}\label{strongly}
\int_{\rr^d}|f(x)|^2dx \leq C\Big(\int_{\rr^d \setminus S}|f(x)|^2dx  
+\int_{\rr^d \setminus \Sigma}|\widehat{f}(\xi)|^2d\xi\Big).
\end{equation*}
\end{itemize}
\end{definition}

\medskip
It can be readily checked that a pair $(S,\Sigma)$ is a strong annihilating one if and only if there exists a positive constant 
$\tilde{C}=\tilde{C}(S,\Sigma)>0$ such that
\begin{equation}\label{strongly2}
\forall f \in L^2(\rr^d), \ \supp\widehat{f}  \subset \Sigma,  \quad \|f\|_{L^2(\rr^d)} \leq \tilde{C}\|f\|_{L^2(\rr^d \setminus S)}.
\end{equation}
We deduce from the Hilbert uniqueness method, Proposition~\ref{ineqspec} and \eqref{cutoff_operator} that the exact controllability of the fractional free Schr\"odinger equation 
\begin{equation*}\label{obs2}
\left\lbrace \begin{array}{ll}
\partial_tf(t,x) + i(-\Delta_x)^{s}f(t,x)=u(t,x)\un_{\omega}(x) \,, \quad &  x \in \mathbb{R}^d,\ t>0, \\
f|_{t=0}=f_0 \in L^2(\rr^d),                                       &  
\end{array}\right.
\end{equation*}
with $s>0$, from a measurable control set $\omega$ holds if and only if its complement $\rr^d \setminus \omega$ forms a strong annihilating pair with the frequency set
$$\mathcal{C}_D(\lambda)= \big\{ \xi \in \rr^d: \ ||\xi|^{2s}-\lambda| \leq \sqrt{D} \big\},$$ 
for some $D>0$ and all $\lambda \in \rr$ with an uniform constant $\tilde{C}$ with respect to $\lambda$ in \eqref{strongly2}.
An exhaustive description of annihilating pairs is for now out of reach. However, the Logvinenko-Sereda Theorem \cite{Logvinenko_Sereda} gives a complete description of all support sets forming a strong annihilating pair with a given bounded frequency set:

\medskip

\begin{theorem}
[Logvinenko-Sereda {\cite{Logvinenko_Sereda}}]
\label{Logvinenko-Sereda}
Let $S,\Sigma\subset\rr^d$ be measurable subsets with $\Sigma$ bounded. Denoting $\tilde S=\rr^d\setminus S$,
the following assertions are equivalent:
\begin{itemize}
\item[-] The pair $(S,\Sigma)$ is a strong annihilating pair
\item[-] The subset $\tilde{S}$ is thick, that is, there exists a cube $K \subset \rr^d$ with sides parallel to coordinate axes and a positive constant $0<\gamma \leq 1$ such that 
$$\forall x \in \rr^d, \quad |(K+x) \cap \tilde{S}| \geq \gamma|K|>0,$$
where $|A|$ denotes the Lebesgue measure of the measurable set $A$.
\end{itemize}
\end{theorem}

\medskip

Notice that if the support set $S$ forms a strong annihilating pair with some bounded frequency set $\Sigma$, then $S$ forms a strong annihilating pair with any bounded frequency set. In order to be able to use the Logvinenko-Sereda theorem in control theory, it is essential to understand how the constant $\tilde{C}(S, \Sigma)$ appearing in \eqref{strongly2} depends on the bounded frequency set $\Sigma$. In \cite[Theorem~3]{Kovrijkine}, Kovrijkine adresses this question by establishing the following quantitative version of the Logvinenko-Sereda Theorem:

\medskip

\begin{theorem}[Kovrijkine {\cite[Theorem~3]{Kovrijkine}}]\label{Kovrijkine1} Let $\omega \subset \rr^d$ be a measurable subset $\gamma$-thick at scale $L>0$, that is, satisfying 
 \begin{equation}\label{thick_def}
 \forall x \in \rr^d, \quad |([0,L]^d+x) \cap \omega| \geq \gamma L^d>0, 
 \end{equation}
with  $0 < \gamma \leq 1$. There exists a universal positive constant $C>0$ independent on the dimension $d \geq 1$ such that for all $f \in L^2(\rr^d)$ satisfying $\supp \hat{f} \subset J$, with $J$ a cube with sides of length $b$ parallel to coordinate axes, 
\begin{equation}\label{kovrijkine1.1}
\|f\|_{L^2(\rr^d)} \leq c(\gamma,d, L, b) \|f\|_{L^2(\omega)},
\end{equation} 
with $$c(\gamma, d, L, b)= \Big( \frac{C^d}{\gamma} \Big)^{Cd(Lb+1)}.$$
\end{theorem}

\medskip

This result is extended to the case of frequency supports covered by the union of finitely many parallepipeds in \cite[Theorem~4]{Kovrijkine}. In this work, we only use this result in the one-dimensional setting: 

\medskip

\begin{theorem}[Kovrijkine {\cite[Theorem~2]{Kovrijkine}}]{\label{Kovrijkine2}} Let $\omega \subset \rr$ be a measurable subset $\gamma$-thick at scale $L>0$, that is, satisfying 
 \begin{equation*}
 \forall x \in \rr, \quad |([0,L]+x) \cap \omega| \geq \gamma L>0,
 \end{equation*}
 with $0 < \gamma \leq 1$.
There exists a universal positive constant $C>0$ such that for all $f \in L^2(\rr)$ satisfying $\supp \hat{f} \subset \bigcup_{k=1}^m J_k$, with $(J_k)_{1 \leq k \leq m}$ a finite family of intervals with length $|J_k|=b$ for all $1 \leq k \leq m$,
\begin{equation}\label{kovrijkine1.2}
\|f\|_{L^2(\rr)} \leq c(\gamma, m, L, b) \|f\|_{L^2(\omega)},
\end{equation} 
with $$c(\gamma, m, L, b)= \Big( \frac{C}{\gamma} \Big)^{Lb \big(\frac{C}{\gamma} \big)^m+m-\frac{1}{2}}.$$
\end{theorem}

\medskip

Notice that the constant $c(\gamma, m, L, b)$ depends only on the number of intervals and their length but not on their locations. Thanks to this explicit dependence of the constant with respect to the length of the intervals in \eqref{kovrijkine1.2}, Egidi and Veseli\'c \cite{veselic}; and Whang, Whang, Zhang and Zhang \cite{Wang} have independently established that the heat equation  
\begin{equation}\label{heat}
\left\lbrace \begin{array}{ll}
(\partial_t -\Delta_x)f(t,x)=u(t,x)\un_{\omega}(x)\,, \quad &  x \in \mathbb{R}^d,\ t>0, \\
f|_{t=0}=f_0 \in L^2(\rr^d),                                       &  
\end{array}\right.
\end{equation}
is null-controllable in any positive time $T>0$ from a measurable control set $\omega \subset \rr^d$ if and only if the control subset $\omega$ is thick in $\rr^d$. In the recent work \cite{kkj}, Beauchard, Jaming and Pravda-Starov prove that this thickness condition is also sufficient for the null-controllability of the harmonic heat equation 
\begin{equation*}\label{harmonic_heat}
\left\lbrace \begin{array}{ll}
(\partial_t +(-\Delta_x +|x|^2))f(t,x)=u(t,x)\un_{\omega}(x)\,, \quad &  x \in \mathbb{R}^d,\ t>0, \\
f|_{t=0}=f_0 \in L^2(\rr^d).                                       &  
\end{array}\right.
\end{equation*}
More generally, the result of \cite[Theorem~2.2]{kkj} shows that the thickness condition is a sufficient condition for the null-controllability of a large class of hypoelliptic quadratic equations.

When the operator $A$ appearing in the evolution equation \eqref{obs01} has a compact resolvent, the result of Proposition~\ref{ineqspec} is specified further as the following infinite dimensional version of the Hautus test:

\medskip

\begin{proposition}[Miller {\cite[Corollary 2.18]{miller}}]\label{hautus_test}
Let $A$ be a selfadjoint operator on $L^2(\rr^d)$, which is the infinitesimal generator of a strongly continuous group $(e^{itA})_{t \in \rr}$ on $L^2(\rr^d)$.
When the operator $A$ has a compact resolvent with a spectral gap $\gamma >0$ in the following sense: $|\lambda - \mu | \geq \gamma$ for all distinct eigenvalues $\lambda$ and $\mu$, the evolution equation \eqref{obs01} is exactly observable at some time $T>0$ from a mesurable subset $\omega$ if and only if 
\begin{equation*}\label{hautus_test2}
\exists \delta>0, \quad \text{for all eigenvector } \psi \text{ of } A, \quad \|\psi \|_{L^2(\rr^d)} \leq \delta \|\psi \|_{L^2(\omega)}.
\end{equation*}
\end{proposition}

\medskip

When the spectral gap condition holds, the observability of the evolution equation \eqref{obs01} then depends only on some properties of the eigenvectors of $A$. Notice in particular that this spectral gap condition holds for the one-dimensional harmonic oscillator $$\mathcal{H}=-\Delta_x+x^2,$$ 
as well as the fractional harmonic oscillators $\mathcal{H}^s$, when $s \geq 1$, since the spectra of these operators are given by 
\begin{equation*}
\sigma\big( (-\Delta_x+x^2)^s \big)= \{(2n+1)^s : \  n \in \nn \}.
\end{equation*}

In this work, we aim at studying necessary and sufficient geometric conditions on the control subsets $\omega$ to ensure exact controllability of the free fractional Schr\"odinger equations
\begin{equation*}\label{free_schrodinger}
\left\lbrace \begin{array}{ll}
\partial_tf(t,x) + i(-\Delta_x)^{s}f(t,x)=u(t,x)\un_{\omega}(x) \,, \quad &  x \in \mathbb{R}^d,\ t>0, \\
f|_{t=0}=f_0 \in L^2(\rr^d),                                       &  
\end{array}\right.
\end{equation*} and the fractional harmonic Schr\"odinger equations
\begin{equation}\label{harmonic_schrodinger}
\left\lbrace \begin{array}{ll}
\partial_tf(t,x) + i(-\Delta_x+|x|^2)^{s}f(t,x)=u(t,x)\un_{\omega}(x) \,, \quad &  x \in \mathbb{R}^d,\ t>0, \\
f|_{t=0}=f_0 \in L^2(\rr^d),                                      &  
\end{array}\right.
\end{equation}
with $s>0$. The first result contained in this paper (Corollary~\ref{gj1}) shows that the thickness condition defined in Theorem~\ref{Logvinenko-Sereda} is a necessary and sufficient geometric condition for the exact controllability of the fractional free Schr\"odinger equations when the fractional parameter satisfies $s \geq \frac{1}{2}$. In the multidimensional setting, the thickness condition is always a necessary condition for exact controllability when the fractional parameter satisfies $s>0$ (Theorem~\ref{prop_necessary_cond}). When the control subset is assumed to be open and $0<s<\frac{1}{2}$, the previous necessary condition is strengthened in Corollary~\ref{gj1rrff} by requiring the complement of the control subset $\rr^d \setminus \omega$ to have empty interior. A sufficient geometric condition to ensure the exact controllability from a neighborhood of the control subset is also given in the multidimensional case when $s \geq \frac{1}{2}$ (Proposition~\ref{gj2}). These results are based on the Logvinenko-Sereda Theorem, its quantitative versions proved by Kovrijkine and new uncertainty principles recently established by Green, Jaye and Mitkovski~\cite{GJM}.

Regarding the fractional harmonic Schr\"odinger equations with $s \geq 1$, we prove that the following geometric condition 
\begin{equation}\label{GC1}
\forall A \in O(\rr^d), \quad \liminf_{R_1 \to +\infty}...\liminf_{R_{d} \to +\infty} \frac{|A(\omega) \cap [-R_1,R_1] \times ... \times [-R_d,R_d]|}{|[-R_1,R_1] \times ... \times [-R_d,R_d]|} >0,
\end{equation}
turns out to be necessary for exact controllability in any dimension $d \geq 1$ (Theorem~\ref{negative_result}) and sufficient in the one-dimensional setting (Theorem~\ref{harmonique}). The proof of this result is based on the infinite dimensional version of the Hautus test given by Proposition~\ref{hautus_test} applied with Hermite functions together with the Rotach-Plancherel formula for Hermite polynomials. This necessary condition enables us to derive a negative result of exact controllability from a cone in dimension~$2$. In the case of the one-dimensional harmonic Schr\"odinger equation, we establish that exact observability from a control subset satisfying \eqref{GC1} holds at any time $T \geq \pi$.

\section{Statements of the main results}

\subsection{Hautus test and exact controllability of fractional harmonic Schr\"odinger equations}

In this section, we state the results related to the exact controllability of fractional harmonic Schr\"odinger equations. We recall that the exact controllability of the evolution equations \eqref{harmonic_schrodinger} is equivalent to the observability of the adjoint systems
\begin{equation}\label{obs_harmonic}
\left\lbrace \begin{array}{ll}
\partial_tf(t,x) - i(-\Delta_x+|x|^2)^{s}f(t,x)=0 \,, \quad &  x \in \mathbb{R}^d,\ t>0, \\
f|_{t=0}=f_0 \in L^2(\rr^d).                                      &  
\end{array}\right.
\end{equation}
The observability of \eqref{obs_harmonic} is deduced thanks to Proposition~\ref{hautus_test} by an infinite dimensional Hautus test on the Hermite functions seen as eigenvectors of the fractional harmonic oscillator.

Hermite functions are normalized $L^2(\rr)$-functions defined as

\begin{equation}\label{e1}
\forall n \geq 0, \quad \psi_n(x)=\frac{1}{\pi^{\frac{1}{4}}2^{\frac{n}{2}} \sqrt{n!}}e^{-\frac{x^2}{2}}H_n(x), 
\end{equation}
where $H_n$ stands for the Hermite polynomial of degree $n \geq 0$, see e.g. \cite{szego}, formula (5.5.1). The family $(\psi_n)_{n \geq 0}$ defines a Hilbert basis of $L^2(\rr)$. They satisfy the following infinite dimensional Hautus test:

\medskip

\begin{proposition}\label{prop}
Let $\omega$ be a measurable subset of $\rr$ and $(\psi_n)_{n \geq 0}$ be the Hermite basis on $L^2(\rr)$ defined in \emph{(\ref{e1})}. The following assertions are equivalent:
\begin{align*}
(i) & \qquad \liminf_{R \to +\infty} \frac{|\omega \cap [-R,R]|}{| [-R,R]|}>0,\\
(ii) & \qquad \exists c>0, \forall n \geq 0, \quad 1=\|\psi_n\|_{L^2(\rr)} \leq c\|\psi_n\|_{L^2(\omega)},\\
(iii) & \qquad \liminf_{n \to +\infty} \|\psi_n\|_{L^2(\omega)}>0.
\end{align*}
\end{proposition}

\medskip

The proof of Proposition~\ref{prop} is given in Section~\ref{proof_prop1} and is based on the Plancherel-Rotach formula.
We readily deduce from Propositions~\ref{hautus_test} and \ref{prop} the following result of exact controllability in the one-dimensional case:

\medskip

\begin{theorem}\label{harmonique}
Let $s \geq 1$ and $\omega$ be a measurable set in $\rr$. The two following assertions are equivalent:
\begin{itemize}
\item[$(i)$] The fractional harmonic Schr\"odinger equation
$$\left\lbrace \begin{array}{ll}
\partial_tf(t,x) + i(-\Delta_x+x^2)^sf(t,x)=u(t,x)\un_{\omega}(x)\,, \quad &  x \in \mathbb{R}, \ t >0, \\
f|_{t=0}=f_0 \in L^2(\rr),                                       &  
\end{array}\right.$$
is exactly controllable from the set $\omega$ for some positive time $T>0$,
\item[$(ii)$] The control set obeys the following geometric condition
$$ \liminf_{R \to +\infty} \frac{|\omega \cap [-R,R]|}{|[-R,R]|}>0.$$
\end{itemize}
If one of the two above conditions holds in the non-fractional case $s=1$, then the harmonic Schr\"odinger equation is exactly controllable from the control subset $\omega$ in any time $T \geq \pi$.
\end{theorem}

\medskip

The equivalence between $(i)$ and $(ii)$ directly follows from Propositions~\ref{hautus_test} and \ref{prop}, since the Hermite basis $(\psi_n)_{n \geq 0}$ is a Hilbert basis of eigenvectors for the fractional harmonic oscillator $\mathcal{H}^s$ associated to the eigenvalues $((2n+1)^s)_{n \geq 0}$ that satisfy the spectral gap condition when $s \geq 1$.

Let us now check in the non-fractional case $s=1$ that when exact controllability holds for some positive time $T_0>0$ then it necessary holds in any time $T \geq \pi$. To that end, we just notice that the function 
$$t\in \rr \mapsto \|e^{it(-\Delta_x+x^2)}f_0\|^2_{L^2(\omega)},$$ 
is $\pi$-periodic since
\begin{multline*}
\forall t \in \rr, \forall f_0 \in L^2(\rr), \quad e^{i(t+\pi)(-\Delta_x+x^2)}f_0 = \sum_{n \in \nn} e^{i(t+\pi) (2n+1)} \langle f_0, \psi_n \rangle_{L^2(\rr)} \psi_n \\
= -e^{it(-\Delta_x+x^2)}f_0.
\end{multline*}
By recalling that exact controllability at some time $T_0>0$ is equivalent to the observability of the evolution equation
 
\begin{equation*}\label{schrodinger2}\left\lbrace \begin{array}{ll}
\partial_tf(t,x) - i(-\Delta_x+x^2)f(t,x)=0 , \quad &  x \in \mathbb{R},\ t>0, \\
f|_{t=0}=f_0 \in L^2(\rr),                                       &  
\end{array}\right.
\end{equation*}
in time $T_0>0$, that is,
\begin{equation*}\label{obs1}
\exists C_{T_0} >0, \forall f_0 \in L^2(\rr), \quad \|f_0 \|^2_{L^2(\rr)} \leq C_{T_0} \int_0^{T_0} \|e^{it(-\Delta_x+x^2)} f_0 \|^2_{L^2(\omega)} dt.
\end{equation*}
The $\pi$-periodicity property readily implies observability and exact controllability at any time $T \geq \pi$, since for all $f_0 \in L^2(\rr)$,
\begin{multline*}
\|f_0 \|^2_{L^2(\rr)} \leq C_{T_0} \int_0^{T_0} \|e^{it(-\Delta_x+x^2)} f_0 \|^2_{L^2(\omega)} dt \leq C_{T_0} \int_0^{(\left\lfloor \frac{T_0}{\pi} \right\rfloor +1) \pi}  \|e^{it(-\Delta_x+x^2)} f_0 \|^2_{L^2(\omega)} dt \\
\leq C_{T_0} \sum_{k=0}^{\lfloor \frac{T_0}{\pi} \rfloor} \int_{k\pi}^{(k +1) \pi}  \|e^{it(-\Delta_x+x^2)} f_0 \|^2_{L^2(\omega)} dt \leq C_{T_0} \Big(\Big\lfloor \frac{T_0}{\pi} \Big\rfloor+1\Big) \int_0^{\pi}  \|e^{it(-\Delta_x+x^2)} f_0 \|^2_{L^2(\omega)} dt,
\end{multline*}
where $\left\lfloor \cdot \right\rfloor$ denotes the floor function.

\medskip
\begin{remark}
When the fractional parameter satisfies $0< s < 1$ and $\omega$ is a measurable set in $\rr$, the geometric condition
$$ \liminf_{R \to +\infty} \frac{|\omega \cap [-R,R]|}{|[-R,R]|}>0,$$
is a necessary condition to ensure the exact controllability of the fractional harmonic Schr\"odinger equation
$$\left\lbrace \begin{array}{ll}
\partial_tf(t,x) + i(-\Delta_x+x^2)^sf(t,x)=u(t,x)\un_{\omega}(x)\,, \quad &  x \in \mathbb{R}, \ t >0, \\
f|_{t=0}=f_0 \in L^2(\rr),                                       &  
\end{array}\right.$$
from the control subset $\omega$ at some positive time $T>0$. Indeed, the a priori estimate characterizing the exact observability property applied to Hermite functions readily implies that condition~$(ii)$ in Proposition~\ref{prop} holds. The same proposition then shows that the above geometric condition has to hold when the fractional harmonic Schr\"odinger equation is exactly controllable.
\end{remark}
\medskip

Theorem~\ref{harmonique} provides a necessary and sufficient geometric condition on the control subset $\omega$ to ensure exact controllability of fractional harmonic Schr\"odinger equations in the one-dimensional setting. In higher dimensions $d \geq 2$, it is natural to investigate if the geometric condition
\begin{equation}\label{wrong_condi}
\liminf_{R \to +\infty} \frac{|\omega \cap B_d(0,R)|}{|B_d(0,R)|}>0,
\end{equation}
where $B_d(0,R)$ denotes the Euclidean ball in $\rr^d$ centered at $0$ with radius $R$, turns out to be sufficient to ensure exact controllability. This is actually not the case in general as the next result shows that condition \eqref{wrong_condi}  is not sufficient for the exact controllability of the harmonic Schr\"odinger equation posed on~$\rr^2$. The stronger necessary geometric condition for exact controllability given by the following proposition proves in particular that the exact controllability of the harmonic Schr\"odinger equation posed on~$\rr^2$
cannot hold from the cone 
\begin{equation}\label{cone}
\omega_{\delta}= \Big\{(r\cos\theta,r\sin\theta): \  r \geq0,\ |\theta| \leq \frac{\pi}{2} -\delta\Big\} \subset \rr^2,
\end{equation}
with $0 < \delta \leq \frac{\pi}{2}$,\\

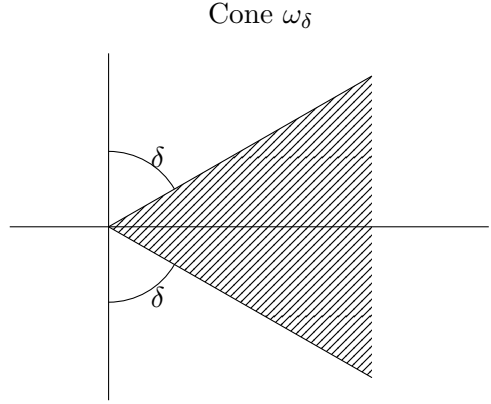
\bigskip

\begin{center}
\begin{tikzpicture}[scale=1]
\draw (2,2.5) node[above]{Cone $\omega_{\delta}$};
\draw (0,0)--({2*sqrt(3)},2);
\draw (0,0)--({2*sqrt(3)},-2);
\draw (0,-2.3)->(0,2.3);
\draw (-1.3,0)->(5,0);
\draw ({0.5*sqrt(3)},0.5) arc (30:90:1);
\draw ({0.5*sqrt(3)},-0.5) arc (-30:-90:1);
\draw (0.65,1.2) node[below]{$\delta$};
\draw (0.65,-1.2) node[above]{$\delta$};
\fill [pattern=north east lines] (0,0)--({2*sqrt(3)},2)--({2*sqrt(3)},-2);
\end{tikzpicture}  
\end{center}

\bigskip

\noindent as the control subset $\omega_{\delta}$ satisfies condition \eqref{wrong_condi}, but fails to satisfy the following necessary geometric condition $$\liminf_{R_1 \to +\infty} \liminf_{R_{2} \to +\infty} \frac{|\omega \cap [-R_1,R_1] \times[-R_2,R_2]|}{|[-R_1,R_1] \times[-R_2,R_2]|} >0,$$ given by the following result:

\medskip

\begin{theorem}\label{negative_result}
Let $s \geq 1$ and $\omega$ be a measurable subset of $\rr^d$ with $d \geq 1$. If the fractional  harmonic Schr\"odinger equation 
\begin{equation}\label{harmonic_control}
\left\lbrace \begin{array}{ll}
\partial_tf(t,x) + i(-\Delta_x+|x|^2)^sf(t,x)=u(t,x)\un_{\omega}(x)\,, \quad &  x \in \mathbb{R}^d,\ t>0, \\ 
f|_{t=0}=f_0 \in L^2(\rr^d),                                       &  
\end{array}\right.
\end{equation}
is exactly controllable from $\omega$ for some positive time $T>0$, then the control subset $\omega$ does obey the following geometric condition
\begin{equation}\label{necessary_cond}
\forall A \in O(\rr^d), \quad \liminf_{R_1 \to +\infty}...\liminf_{R_d \to +\infty} \frac{|A(\omega) \cap [-R_1,R_1] \times ... \times [-R_d,R_d]|}{|[-R_1,R_1] \times ... \times [-R_d,R_d]|} >0,
\end{equation}
where $O(\rr^d)$ stands for the orthogonal group in dimension $d$.
\end{theorem}

\medskip

This result shows that the exact controllability of the fractional harmonic Schr\"odinger equations requires the control subset to be distributed in any space direction. The sufficiency of condition \eqref{necessary_cond} for exact controllability is an open question in dimension $d \geq 2$. The proof of Theorem~\ref{negative_result} is given in Section~\ref{proof_negative_result}. An example of sufficient geometric condition is given by the following proposition, which is deduced from Theorem~\ref{harmonique}:

\medskip

\begin{exemple}
Let $d \geq 1$, $s \geq 1$ and $\omega$ be a measurable subset of $\rr$ satisfying 
$$ \liminf_{R \to +\infty} \frac{|\omega \cap [-R,R]|}{|[-R,R]|}>0.$$
The evolution equation
\begin{equation}\label{harmonic_multid}
\left\lbrace \begin{array}{ll}
\partial_tf(t,x) + i(-\Delta_x+|x|^2)^sf(t,x)=u(t,x)\un_{\omega \times \rr^{d-1}}(x)\,, \quad &  x \in \mathbb{R}^d, \ t>0, \\
f|_{t=0}=f_0 \in L^2(\rr^d),                                       &  
\end{array}\right.
\end{equation}
is exactly controllable from $\omega \times \rr^{d-1}$ for some positive time $T>0$.
\end{exemple}

\medskip

\begin{proof}
We first notice that the spectral gap condition in Proposition~\ref{hautus_test} is satisfied by the fractional harmonic operator $(-\Delta_x + |x|^2)^s$, since its spectrum is given by 
$$\sigma\big((-\Delta_x +|x|^2)^s\big) =\big\{(2n+d)^s:\ n \in \nn\big\}.$$
Let $N \in \nn$ and $f$ be an eigenvector associated to the eigenvalue $(2N+d)^s$. There exists a family of complex numbers $(C_{\alpha})_{\alpha \in \nn^d,  |\alpha|=N}$ such that 
$$f = \sum \limits_{\substack{\alpha \in \nn^d, \\ |\alpha|=N}} C_{\alpha} \Psi_{\alpha},$$ 
with 
\begin{multline}\label{bbb2}
\forall \alpha=(\alpha_1,...,\alpha_d)=(\alpha_1,\alpha') \in \nn^d, \forall x=(x_1,...,x_d)=(x_1,x') \in \rr^d,\\  \Psi_{\alpha}(x)= \prod_{i=1}^d \psi_{\alpha_i} (x_i).
\end{multline} 
We notice that
\begin{align*}
& \ \|f\|_{L^2(\omega \times \rr^{d-1})}^2 = \int_{\omega} \Big(\int_{\rr^{d-1}} \Big| \sum \limits_{\substack{\alpha \in \nn^d,\\ |\alpha|=N}} C_{\alpha} \Psi_{\alpha}(x_1,x')\Big|^2 dx'\Big)dx_1 \\
= & \ \int_{\omega} \Big(\int_{\rr^{d-1}} \Big| \sum \limits_{\substack{\alpha' \in \nn^{d-1},\\ |\alpha'| \leq N}} C_{N-|\alpha'|, \alpha'} \psi_{N-|\alpha'|}(x_1) \Psi_{\alpha'}(x')\Big|^2 dx'\Big)dx_1 \\
= & \ \int_{\omega} \sum \limits_{\substack{\alpha' \in \nn^{d-1},\\ |\alpha'| \leq N}} |C_{N-|\alpha'|, \alpha'}|^2 |\psi_{N-|\alpha'|}(x_1)|^2 dx_1 = \sum \limits_{\substack{\alpha' \in \nn^{d-1},\\ |\alpha'| \leq N}} |C_{N-|\alpha'|, \alpha'} |^2 \|\psi_{N-|\alpha'|} \|_{L^2(\omega)}^2.
\end{align*}
By using Proposition~\ref{prop},
\begin{equation*}
\exists c >0, \forall n \geq 0, \quad 1= \|\psi_n\|_{L^2(\rr)} \leq c \| \psi_n \|_{L^2(\omega)},
\end{equation*}
we deduce that
\begin{equation}\label{ht}
\|f\|_{L^2(\omega \times \rr^{d-1})}^2 \geq \frac{1}{c^2} \sum \limits_{\substack{\alpha' \in \nn^{d-1},\\ |\alpha'| \leq N}}|C_{N-|\alpha'|, \alpha'}|^2=
\frac{1}{c^2} \sum \limits_{\substack{\alpha \in \nn^{d},\\ |\alpha|= N}}|C_{\alpha}|^2
= \frac{1}{c^2} \| f\|_{L^2(\rr^d)}^2.
\end{equation}
Thanks to Proposition~\ref{hautus_test}, we deduce from \eqref{ht} the observability of the adjoint system 
\begin{equation*}
\left\lbrace \begin{array}{ll}
\partial_tg(t,x) - i(-\Delta_x+|x|^2)^sg(t,x)=0\,, \quad &  x \in \mathbb{R}^d, \ t>0, \\
g|_{t=0}=g_0 \in L^2(\rr^d),                                       &  
\end{array}\right.
\end{equation*}
and then the exact controllability of the evolution equation \eqref{harmonic_multid}.
\end{proof}

The above example ensures the exact controllability of fractional harmonic Schr\"odinger equations from the half space $\rr_+ \times \rr^{d-1}$. In particular, it appears that the cone (\ref{cone}) is a limit case for exact controllability.

\subsection{Exact controllability of fractional free Schr\"odinger equations}
In this section, we study the exact controllability of the fractional free Schr\"odinger equation
\begin{equation}\label{schrodinger}
\left\lbrace \begin{array}{ll}
(\partial_t +i(-\Delta_x)^s)f(t,x)=u(t,x)\un_{\omega}(x)\,, \quad &  x \in \mathbb{R}^d,\ t>0, \\
f|_{t=0}=f_0 \in L^2(\rr^d),                                       &  
\end{array}\right.
\end{equation}
with $s >0$, from a measurable control subset $\omega \subset \rr^d$.
By the Hilbert uniqueness method, the exact controllability of the evolution equation (\ref{schrodinger}) is equivalent to the observability of the adjoint system
\begin{equation}\label{schrodinger*}
\left\lbrace \begin{array}{ll}
(\partial_t -i(-\Delta_x)^s)g(t,x)=0, \quad &  x \in \mathbb{R}^d,\ t>0, \\
g|_{t=0}=g_0 \in L^2(\rr^d).                                     &  
\end{array}\right.
\end{equation}
We first notice that the Logvinenko-Sereda Theorem recalled in Theorem~\ref{Logvinenko-Sereda} shows that the thickness condition is a necessary geometric condition for the exact controllability of the fractional Schr\"odinger equation~\eqref{schrodinger} in any dimension $d \geq 1$: 

\medskip

\begin{theorem}\label{prop_necessary_cond} Let $s >0$, $\omega$ be a measurable subset of $\rr^d$ and $T>0$. If the fractional free Schr\"odinger equation \emph{(\ref{schrodinger})} is exactly controllable from $\omega$ at some time $T>0$, then the control subset $\omega$ is a thick set.
\end{theorem}

\medskip

\begin{proof}
The assumptions of Theorem~\ref{prop_necessary_cond} imply that the system (\ref{schrodinger*}) is exactly observable from $\omega$ in time $T>0$. We deduce from Proposition~\ref{ineqspec} that there exist some positive constants $k>0$ and $D>0$ such that for all $\lambda \in \rr$ and $f \in L^2(\rr^d)$,
\begin{equation}\label{spec1}
\supp\widehat{f} \subset \big\{ \xi \in \rr^d:\ ||\xi|^{2s}- \lambda | \leq \sqrt{D} \big\} \implies \|f\|_{L^2(\rr^d)} \leq \sqrt{k} \|f\|_{L^2(\omega)}.
\end{equation}
While taking $\lambda=0$, the assertion \eqref{spec1} implies in particular that the sets $\rr^d \setminus \omega$ and the Euclidean ball $B_d(0, D^{\frac{1}{4s}})$ in $\rr^d$ make a strong annihilating pair. Since $B_d(0, D^{\frac{1}{4s}})$ is a bounded set, the Logvinenko-Sereda Theorem (Theorem~\ref{Logvinenko-Sereda}) ensures that the control subset $\omega$ must be thick.
\end{proof}

In the one-dimensional setting, the quantitative version of the Logvinenko-Sereda Theorem established by Kovrijkine (Theorem~\ref{Kovrijkine1}) together with Proposition~\ref{ineqspec} show that the thickness condition is also sufficient for the exact controllability of the fractional free Schr\"odinger equations (\ref{schrodinger}) when $s \geq \frac{1}{2}$: 

\medskip

\begin{theorem}\label{thick_control}
If $s \geq \frac{1}{2}$ and $\omega$ is a thick measurable subset of $\rr$, then there exists a positive constant $T_0 >0$ such that for any time $T>T_0$, the fractional free Schr\"odinger equation \emph{(\ref{schrodinger})} is exactly controllable from $\omega$ in time $T$.
\end{theorem}

\medskip

\begin{proof}
Let $\omega \subset \rr$ be a measurable subset $\gamma$-thick at scale $L>0$.
Let $D >0$, $\lambda \in \rr$ and $f \in L^2(\rr)$ such that 
$$\supp \widehat{f} \subset \big\{\xi \in \rr: \ ||\xi|^{2s}-\lambda| \leq \sqrt{D}\big\}.$$ 
If $|\lambda| > \sqrt{D}+1$, the Fourier transform $\widehat{f}$ is supported in the union of two intervals 
\begin{equation*}
\supp \hat{f} \subset \big[-(\sqrt{D}+|\lambda|)^{\frac{1}{2s}}, -(-\sqrt{D}+|\lambda|)^{\frac{1}{2s}}\big] \cup \big[(-\sqrt{D}+|\lambda|)^{\frac{1}{2s}}, (\sqrt{D}+|\lambda|)^{\frac{1}{2s}}\big]
\end{equation*}
and 
\begin{equation*}
\big|\big[-(\sqrt{D}+|\lambda|)^{\frac{1}{2s}}, -(-\sqrt{D}+|\lambda|)^{\frac{1}{2s}}\big] \big|=\big|\big[(-\sqrt{D}+|\lambda|)^{\frac{1}{2s}}, (\sqrt{D}+|\lambda|)^{\frac{1}{2s}}\big] \big| \leq \frac{\sqrt{D}}{s},
\end{equation*}
since $s \geq \frac{1}{2}$ and $|\lambda| \geq \sqrt{D}+1$. It follows from Theorem~\ref{Kovrijkine2} that there exists a universal positive constant $C'>0$ such that 
\begin{equation}\label{kov1}
\|f\|_{L^2(\rr)} \leq \Big( \frac{C'}{\gamma} \Big)^{L\frac{\sqrt{D}}{s}(\frac{C'}{\gamma})^2+\frac{3}{2}} \|f\|_{L^2(\omega)}.
\end{equation} 
On the other hand, we notice that $\supp \hat{f} \subset \big[0, \big(2\sqrt{D}+1 \big)^{\frac{1}{2s}} \big]$ when $|\lambda| \leq \sqrt{D}+1$. We deduce from Theorem~\ref{Kovrijkine1} that there exists a universal positive constant $C >0$ such that 

\begin{equation}\label{kov2}
\|f\|_{L^2(\rr)} \leq \Big(\frac{C}{\gamma}\Big)^{C(1+L(2\sqrt{D}+1)^{\frac{1}{2s}})}\|f\|_{L^2(\omega)}.
\end{equation}
By setting 
$$d= \max \Big(\Big( \frac{C'}{\gamma} \Big)^{L\frac{\sqrt{D}}{s}(\frac{C'}{\gamma})^2+\frac{3}{2}}, \Big(\frac{C}{\gamma}\Big)^{C(1+L(2\sqrt{D}+1)^{\frac{1}{2s}})}\Big)^2,$$ 
it follows from \eqref{kov1} and \eqref{kov2} that for all $\lambda \in \rr$ and $f \in L^2(\rr)$ with $\supp \hat{f} \subset \big\{ \xi \in \rr:\ ||\xi|^{2s}-\lambda| \leq \sqrt{D}\big\}$, 
\begin{equation*}
\| f\|_{L^2(\rr)} \leq \sqrt{d} \|f\|_{L^2(\omega)}.
\end{equation*}
Theorem~\ref{thick_control} is then a direct consequence of the Hilbert uniqueness method and Proposition~\ref{ineqspec}.
\end{proof}

We deduce the following necessary and sufficient geometric condition for the exact controllability of the fractional free Schr\"odinger equations:

\medskip

\begin{corollary}\label{gj1}
Let $s \geq \frac{1}{2}$ and $\omega$ be a measurable subset of $\rr$. The one-dimensional fractional free Schr\"odinger equation 
\begin{equation*}
\left\lbrace \begin{array}{ll}
(\partial_t +i(-\Delta_x)^s)f(t,x)=u(t,x)\un_{\omega}(x)\,, \quad &  x \in \mathbb{R},\ t>0, \\
f|_{t=0}=f_0 \in L^2(\rr),                                       &  
\end{array}\right.
\end{equation*}
is exactly controllable from $\omega$ if and only if $\omega$ is thick.
\end{corollary}

\medskip

It is interesting to point out that the sharp geometric condition for the exact controllability of the free Schr\"odinger equation is actually the same as the one for the null-controllability of the free heat equation \eqref{heat}, even if the behavior of these two equations are very different. This thickness condition is also a necessary and sufficient condition for the null-controllability of the fractional heat equations as showed by Alphonse and Bernier~\cite[Remark~1.13]{alphonse}, when the fractional parameter satisfies $s > \frac{1}{2}$.

By elaborating further on the link made by Duyckaerts and Miller~\cite{duyckaerts} (Corollary~2) between results of exact controllability at some positive time for Schr\"odinger equations and results of null-controllability in any positive time for the associated heat equations, 
the recent result of non-null-controllability in any positive time established by Koenig~\cite{koenig} (Theorem 3) for the fractional heat equation 
\begin{equation}\label{heatfrac}
\left\lbrace \begin{array}{ll}
(\partial_t +(-\Delta_x)^s)f(t,x)=u(t,x)\un_{\omega}(x)\,, \quad &  x \in \mathbb{R}^d,\ t>0, \\
f|_{t=0}=f_0 \in L^2(\rr^d),                                       &  
\end{array}\right.
\end{equation}
when $0<s<\frac{1}{2}$, once $\omega \subset \rr^d$ is an open subset whose complement $\rr^d \setminus \omega$ has a non-empty interior, readily implies the following necessary condition for the exact controllability of fractional free Schr\"odinger equations from an open control subset when the fractional parameter satisfies $0<s<\frac{1}{2}$:

\medskip

\begin{corollary}\label{gj1rrff}
Let $0<s<\frac{1}{2}$ and $\omega \subset \rr^d$ be an open subset. If the fractional free Schr\"odinger equation 
\begin{equation*}
\left\lbrace \begin{array}{ll}
(\partial_t +i(-\Delta_x)^s)f(t,x)=u(t,x)\un_{\omega}(x)\,, \quad &  x \in \mathbb{R}^d,\ t>0, \\
f|_{t=0}=f_0 \in L^2(\rr^d),                                       &  
\end{array}\right.
\end{equation*}
is exactly controllable from $\omega$ at some positive time $T>0$, then the complement of the control subset $\rr^d \setminus \omega$ has empty interior.
\end{corollary}

\medskip

In a recent work~\cite{GJM} , Green, Jaye and Mitkovski prove new uncertainty principles. These results can be seen as a generalization of the Logvinenko-Sereda Theorem and its quantitative versions given by Kovrijkine. They establish in particular the following uncertainty principle:

\medskip

\begin{theorem}[Green, Jaye \& Mitkovski {\cite[Corollary~3]{GJM}}]
Let $\omega$ be a measurable subset of $\rr^d$ satisfying the following one-dimensional geometric control condition: there exist some positive constants $L>0$ and $0<\gamma \leq 1$ such that for all straight lines $\mathcal{D} \subset \rr^d$ and for all line segment $S \subset \mathcal{D}$ of length $L$, 
\begin{equation*}
\mathcal{H}(S \cap \omega) \geq \gamma L,
\end{equation*}
where $\mathcal{H}$ denotes the one-dimensional Lebesgue measure of $\mathcal{D}$.
Then, for all $\delta >0$, $\beta >0$, there exists $C>0$ such that for all $R>0$ and $f \in L^2(\rr^d)$,
\begin{equation*}
\supp \hat{f} \subset \big\{ \xi \in \rr^d:\ R- \beta \leq |\xi| \leq R+ \beta \big\} \implies \|f \|_{L^2(\rr^d)} \leq C \|f \|_{L^2(\mathcal{U}_{\delta}(\omega))},
\end{equation*}
where $\mathcal{U}_{\delta} (\omega)$ denotes the $\delta$-neighborhood of $\omega$ defined by $\mathcal{U}_{\delta}(\omega)= \{x \in \rr^d: \ d(x,\omega) < \delta \}$.
\end{theorem}

\medskip

This result together with Proposition~\ref{ineqspec} and the Hilbert uniqueness method allow to derive the following sufficient geometric condition on the control subset $\omega$ to ensure the exact controllability of the fractional free Schr\"odinger equation \eqref{schrodinger}: 

\medskip

\begin{proposition}\label{gj2}
Let $d \geq 1$, $s \geq \frac{1}{2}$ and $\omega$ be a measurable subset of $\rr^d$ satisfying that there exist some positive constants $L>0$ and $0<\gamma \leq 1$ such that for all straight lines $\mathcal{D} \subset \rr^d$ and for all line segment $S \subset \mathcal{D}$ of length $L$, 
\begin{equation}\label{thick1d}
\mathcal{H}(S \cap \omega) \geq \gamma L,
\end{equation}
where $\mathcal{H}$ denotes the one-dimensional Lebesgue measure of $\mathcal{D}$. Then, for all $\delta >0$, the fractional free Schr\"odinger equation 
\begin{equation*}
\left\lbrace \begin{array}{ll}
(\partial_t +i(-\Delta_x)^s)f(t,x)=u(t,x)\un_{\mathcal{U}_{\delta}(\omega)}(x)\,, \quad &  x \in \mathbb{R}^d,\ t>0, \\
f|_{t=0}=f_0 \in L^2(\rr^d),                                       &  
\end{array}\right.
\end{equation*}
is exactly controllable from the control subset $\mathcal{U}_{\delta}(\omega)$.
\end{proposition}

\medskip

Notice that any measurable subset $\omega \subset \rr^d$ satisfying the condition \eqref{thick1d} is a thick subset in the sense of \eqref{thick_def}. The condition \eqref{thick1d} is therefore stronger than the thickness condition. 

\bigskip

\textbf{Addendum.} After the completion of this work, we find out the recent preprint~\cite{huang}, which points out the same necessary and sufficient geometric condition on control subsets for the exact controllability of the one-dimensional free and harmonic Schr\"odinger equations. The authors of~\cite{huang} also derive a necessary and sufficient geometric condition for the exact controllability of the Schr\"odinger equation associated to the one-dimensional operator $P=-\Delta_x+V(x)$, with $V(x)=x^{2m}$, when $m \geq 2$, which is not discussed in this work. On the other hand, the fractional case is not discussed in~\cite{huang} and the multidimensional results contained in this work are also not covered by the results of~\cite{huang}. The proofs given in the two papers are also different. Regarding the free Schr\"odinger equation, the proof of~\cite{huang} given only in the non-fractional case $s=1$ makes the use of explicit computations of the Fourier transform of Gaussians which cannot be directly used in the fractional case. This proof does not make explicit the condition $s \geq \frac{1}{2}$ found out in Corollary~\ref{gj1} to obtain the necessary and sufficient geometric condition for exact controllability. Regarding the harmonic Schr\"odinger equation, the proof of~\cite{huang} is much more involved since it is given in the general case when $m \geq 1$ with no explicit simplification in the harmonic case $m=1$. The use of the Plancherel-Rotach formula in this work allows to obtain more directly the result in the harmonic case and accounts for the fact that the present paper is much shorter than the preprint~\cite{huang}.

\section{Proof of Proposition~\ref{prop}}\label{proof_prop1}

This section is devoted to the proof of Proposition~\ref{prop}.
The assertions $(ii)$ and $(iii)$ are clearly equivalent. Indeed, $(ii)$ implies $(iii)$ by passing to the limit inferior
$$\liminf_{n \to +\infty} \|\psi_n\|_{L^2(\omega)} \geq \frac{1}{c}>0.$$
On the other hand, we readily deduce from formula (\ref{e1}) that for any measurable set $\omega \subset \rr$ of positive measure $|\omega| >0$,
\begin{equation}\label{rr1}
\forall n \geq 0, \quad \|\psi_n\|_{L^2(\omega)}>0, 
\end{equation}
since $H_n$ is a polynomial of degree $n$. Then, assertion $(ii)$ directly follows from $(iii)$ and (\ref{rr1}). It is therefore sufficient to prove the equivalence of assertions $(i)$ and $(iii)$.

As a preliminary step, we begin by establishing that 
\begin{equation}\label{rr2}
\forall 0<\eps <\frac{\pi}{2}, \quad \lim_{n \to +\infty}\|\psi_n\|_{L^2(I_{\eps,n})}=\sqrt{1-\frac{2\eps}{\pi}},
\end{equation}
where $I_{\eps,n}$ stands for the open interval $(-\sqrt{2n+1}\cos \eps,\sqrt{2n+1}\cos \eps)$.
Let $0<\eps<\frac{\pi}{2}$. The Plancherel-Rotach formula for Hermite polynomials~\cite{szego} (Theorem~8.22.9) provides that for all $\eps \leq \theta \leq \pi-\eps$,
\begin{equation}\label{e2}
e^{-\frac{x^2}{2}}H_n(x)=\frac{2^{\frac{n}{2}+\frac{1}{4}}\sqrt{n!}}{(\pi n)^{\frac{1}{4}}\sqrt{\sin \theta}}\Big[\sin\Big(\Big(\frac{n}{2}+\frac{1}{4}\Big)\big(\sin(2\theta)-2\theta\big)+\frac{3\pi}{4}\Big)+R_n(x)\Big],
\end{equation}
with $x=\sqrt{2n+1}\cos \theta$, where the remainder term satisfies
\begin{equation}\label{e3}
\exists C>0, \forall n \geq 0, \forall  y \in [-\sqrt{2n+1}\cos \eps,\sqrt{2n+1}\cos \eps], \quad |R_n(y)| \leq \frac{C}{n+1}.
\end{equation}
It follows from (\ref{e1}) and (\ref{e2}) that the Hermite function 
\begin{equation*}\label{e30}
\psi_n(x)=\frac{1}{\pi^{\frac{1}{4}}2^{\frac{n}{2}} \sqrt{n!}}e^{-\frac{x^2}{2}}H_n(x), \quad \|\psi_n\|_{L^2(\rr)}=1, \quad n \geq 0,
\end{equation*}
satisfies that for all $\eps \leq \theta \leq \pi-\eps$ and $n\geq 1$,
\begin{equation}\label{e4}
\psi_n(x)=\frac{2^{\frac{1}{4}}}{\sqrt{\pi} n^{\frac{1}{4}}\sqrt{\sin \theta}}\Big[\sin\Big(\Big(\frac{n}{2}+\frac{1}{4}\Big)\big(\sin(2\theta)-2\theta\big)+\frac{3\pi}{4}\Big)+R_n(x)\Big],
\end{equation}
with $x=\sqrt{2n+1}\cos \theta$.
We obtain from (\ref{e4}) that for all $n \geq 1$,
\begin{equation}\label{e6}
\|F_n\|_{L^2(I_{\eps,n})}- \|G_n\|_{L^2(I_{\eps,n})} \leq \|\psi_n\|_{L^2(I_{\eps,n})} 
\leq \|F_n\|_{L^2(I_{\eps,n})}+ \|G_n\|_{L^2(I_{\eps,n})} ,
\end{equation}
where the two functions $F_n$ and $G_n$ are defined for all $\eps < \theta < \pi-\eps$,
\begin{equation}\label{e7}
F_n(x)=\frac{2^{\frac{1}{4}}}{\sqrt{\pi} n^{\frac{1}{4}}\sqrt{\sin \theta}}\sin\Big(\Big(\frac{n}{2}+\frac{1}{4}\Big)\big(\sin(2\theta)-2\theta\big)+\frac{3\pi}{4}\Big)
\end{equation}
and
\begin{equation}\label{e8}
G_n(x)=\frac{2^{\frac{1}{4}}}{\sqrt{\pi} n^{\frac{1}{4}}\sqrt{\sin \theta}}R_n(x),
\end{equation}
with $x=\sqrt{2n+1}\cos \theta$. We deduce from (\ref{e3}) that for all $n \geq 1$,
$$\|G_n\|_{L^2(I_{\eps,n})}^2=\frac{\sqrt{2(2n+1)}}{\pi \sqrt{n}}\int_{\eps}^{\pi-\eps}|R_n(\sqrt{2n+1}\cos \theta)|^2d\theta \leq \frac{C^2\sqrt{2(2n+1)}(\pi-2\eps)}{\pi \sqrt{n}(n+1)^2},$$
implying that 
\begin{equation}\label{e9}
\lim_{n \to +\infty}\|G_n\|_{L^2(I_{\eps,n})}=0.
\end{equation}
On the other hand, we deduce from (\ref{e7}) that for all $n \geq 1$,
\begin{equation}\label{e11}
|F_n(x)|^2=\frac{1}{\sqrt{2n}\pi \sin \theta}\Big[1-\cos\Big(\Big(n+\frac{1}{2}\Big)\big(\sin(2\theta)-2\theta\big)+\frac{3\pi}{2}\Big)\Big].
\end{equation}
It follows from (\ref{e11}) that for all $n \geq 1$,
\begin{multline}\label{e12}
\int_{I_{\eps,n}}|F_n(x)|^2dx = \int_{I_{\eps,n}}\frac{\sqrt{2n+1}}{\sqrt{2n}\pi \sqrt{(2n+1)-x^2}}dx\\
-\frac{\sqrt{2n+1}}{\sqrt{2n}\pi }\int_{\eps}^{\pi-\eps}\cos\Big(\Big(n+\frac{1}{2}\Big)\big(\sin(2\theta)-2\theta\big)+\frac{3\pi}{2}\Big)d\theta,
\end{multline}
since $\sin^2 \theta=1-\frac{x^2}{2n+1}$, when $x=\sqrt{2n+1}\cos \theta$.
We have for all $n \geq 1$,
\begin{multline}\label{e13}
\int_{I_{\eps,n}}\frac{\sqrt{2n+1}}{\sqrt{2n}\pi \sqrt{(2n+1)-x^2}}dx =\frac{1}{\pi}\sqrt{\frac{2n+1}{2n}}\int_{-\cos \eps}^{\cos \eps}\frac{dy}{\sqrt{1-y^2}}\\=\frac{2}{\pi}\sqrt{\frac{2n+1}{2n}}\textrm{arcsin}(\cos \eps)
=\frac{2}{\pi}\sqrt{\frac{2n+1}{2n}}\Big(\frac{\pi}{2}-\textrm{arccos}(\cos \eps)\Big)=\sqrt{\frac{2n+1}{2n}}\Big(1-\frac{2\eps}{\pi}\Big),
\end{multline}
since
$$\forall -1 \leq x \leq 1, \quad \textrm{arccos} \ x+\textrm{arcsin} \ x=\frac{\pi}{2}.$$
It follows from (\ref{e13}) that
\begin{equation}\label{e14}
\lim_{n \to +\infty}\int_{I_{\eps,n}}\frac{\sqrt{2n+1}}{\sqrt{2n}\pi \sqrt{(2n+1)-x^2}}dx=1-\frac{2\eps}{\pi}.
\end{equation}
On the other hand, we have for all $n \geq 1$,
\begin{multline}\label{e15}
\frac{\sqrt{2n+1}}{\sqrt{2n}\pi }\int_{\eps}^{\pi-\eps}\cos\Big(\Big(n+\frac{1}{2}\Big)\big(\sin(2\theta)-2\theta\big)+\frac{3\pi}{2}\Big)d\theta\\
=\frac{\sqrt{2n+1}}{\sqrt{2n}\pi }\int_{\varphi(\eps)}^{\varphi(\pi-\eps)}\frac{1}{\varphi'(\varphi^{-1}(t))}\cos\Big(\Big(n+\frac{1}{2}\Big)t-\frac{3\pi}{2}\Big)dt,
\end{multline}
where $\varphi(\theta)=2\theta-\sin(2\theta)$ is a increasing $C^{\infty}$-diffeomorphism from $[\eps,\pi-\eps]$ to 
$$[\varphi(\eps),\varphi(\pi-\eps)]=[2\eps-\sin(2\eps),2\pi-2\eps+\sin(2\eps)].$$
An integration by parts shows that for all $n \geq 1$,
\begin{multline*}\label{e27}
\int_{\varphi(\eps)}^{\varphi(\pi-\eps)}\frac{1}{\varphi'(\varphi^{-1}(t))}\cos\Big(\Big(n+\frac{1}{2}\Big)t-\frac{3\pi}{2}\Big)dt=\left[\frac{\sin((n+\frac{1}{2})t-\frac{3\pi}{2})}{\varphi'(\varphi^{-1}(t))(n+\frac{1}{2})}\right]_{\varphi(\eps)}^{\varphi(\pi-\eps)}\\
-\frac{1}{n+\frac{1}{2}}\int_{\varphi(\eps)}^{\varphi(\pi-\eps)}\frac{d}{dt}\Big(\frac{1}{\varphi'(\varphi^{-1}(t))}\Big)\sin\Big(\Big(n+\frac{1}{2}\Big)t-\frac{3\pi}{2}\Big)dt,
\end{multline*}
implying that
\begin{equation}\label{e28}
\Big|\int_{\varphi(\eps)}^{\varphi(\pi-\eps)}\frac{1}{\varphi'(\varphi^{-1}(t))}\cos\Big(\Big(n+\frac{1}{2}\Big)t-\frac{3\pi}{2}\Big)dt\Big|
\leq  \frac{C_0}{n+\frac{1}{2}}(2+2\pi-4\eps+2\sin(2\eps)),
\end{equation}
with 
$$0 < C_0=\sup_{t \in [2\eps-\sin(2\eps),2\pi-2\eps+\sin(2\eps)]}\Big\{\Big|\frac{1}{\varphi'(\varphi^{-1}(t))}\Big|+\Big|\frac{d}{dt}\Big(\frac{1}{\varphi'(\varphi^{-1}(t))}\Big)\Big|\Big\}<+\infty.$$
We deduce from (\ref{e15}) and (\ref{e28}) that 
\begin{equation}\label{e29}
\lim_{n \to +\infty}\frac{\sqrt{2n+1}}{\sqrt{2n}\pi }\int_{\eps}^{\pi-\eps}\cos\Big(\Big(n+\frac{1}{2}\Big)\big(\sin(2\theta)-2\theta\big)+\frac{3\pi}{2}\Big)d\theta=0.
\end{equation}
It follows from (\ref{e12}), (\ref{e14}) and (\ref{e29}) that  
\begin{equation}\label{rr4}
\lim_{n \to +\infty}\|F_n\|_{L^2(I_{\eps,n})}=\sqrt{1-\frac{2\eps}{\pi}}.
\end{equation}
Thanks to (\ref{e6}), (\ref{e9}) and (\ref{rr4}), we obtain that (\ref{rr2}) holds. 

Let $\omega$ be a measurable subset of $\rr$ verifying 
\begin{equation}\label{rr5}
\delta=\liminf_{n \to +\infty} \|\psi_n\|_{L^2(\omega)}>0.
\end{equation}
Let $0<\eps_0 <\frac{\pi}{2}$ satisfying
\begin{equation*}\label{rr6}
0<\eps_0 <\frac{\pi}{2}\delta^2.
\end{equation*}
By using that $\|\psi_n\|_{L^2(\rr)}=1$, we observe that for all $0<\eps <\frac{\pi}{2}$ and $n \geq 0$,
\begin{multline}\label{rr7}
\|\psi_n\|_{L^2(\omega)}^2=\|\psi_n\|_{L^2(\omega \cap I_{\eps,n})}^2+\|\psi_n\|_{L^2(\omega \cap (\rr \setminus I_{\eps,n}))}^2 \\ \leq \|\psi_n\|_{L^2(\omega \cap I_{\eps,n})}^2+\|\psi_n\|_{L^2(\rr \setminus I_{\eps,n})}^2
= \|\psi_n\|_{L^2(\omega \cap I_{\eps,n})}^2+1-\|\psi_n\|_{L^2(I_{\eps,n})}^2.
\end{multline}
It follows from (\ref{rr2}), (\ref{rr5}) and (\ref{rr7}) that for all $0<\eps \leq \eps_0$,
\begin{equation}\label{rr8}
\liminf_{n \to +\infty}\|\psi_n\|_{L^2(\omega \cap I_{\eps,n})}>0.
\end{equation}
On the other hand, we deduce from (\ref{e4}), (\ref{e7}) and (\ref{e8}) that for all $n \geq 1$ and $0<\eps <\frac{\pi}{2}$,
\begin{equation}\label{rr9}
\|F_n\|_{L^2(\omega_{\eps,n})}- \|G_n\|_{L^2(\omega_{\eps,n})} \leq \|\psi_n\|_{L^2(\omega_{\eps,n})} 
\leq \|F_n\|_{L^2(\omega_{\eps,n})}+ \|G_n\|_{L^2(\omega_{\eps,n})} ,
\end{equation}
with 
\begin{equation}\label{rr9b}
\omega_{\eps,n}=\omega \cap I_{\eps,n}. 
\end{equation}
By using that $\|G_n\|_{L^2(\omega_{\eps,n})} \leq \|G_n\|_{L^2(I_{\eps,n})}$, it follows from (\ref{e9}) that 
\begin{equation}\label{rr10}
\forall 0<\eps <\frac{\pi}{2}, \quad \lim_{n \to +\infty}\|G_n\|_{L^2(\omega_{\eps,n})}=0.
\end{equation}
We deduce from (\ref{e11}) and (\ref{rr9b}) that for all $n \geq 1$,
\begin{multline}\label{e12b}
\int_{\omega_{\eps_0,n}}|F_n(x)|^2dx \leq \int_{\omega_{\eps_0,n}}\frac{2\sqrt{2n+1}}{\sqrt{2n}\pi \sqrt{(2n+1)-x^2}}dx \\
\leq \int_{\omega_{\eps_0,n}}\frac{\sqrt{2}}{\sqrt{n}\pi \sqrt{1-\cos^2 \eps_0}}dx \leq \frac{\sqrt{2}}{\pi \sin \eps_0 }
\frac{|\omega \cap (-\sqrt{2n+1},\sqrt{2n+1})|}{\sqrt{n}}.
\end{multline}
By using \eqref{rr8}, \eqref{rr9}, \eqref{rr9b}, \eqref{rr10} and \eqref{e12b}, it follows that
\begin{multline}\label{rr11}
\liminf_{n \to +\infty} \frac{|\omega \cap (-\sqrt{2n+1},\sqrt{2n+1})|}{\sqrt{n}} \geq \liminf_{n \to +\infty} \frac{\pi \sin \eps_0}{\sqrt{2}} \int_{\omega_{\eps_0,n}}|F_n(x)|^2dx \\
= \liminf_{n \to +\infty} \frac{\pi \sin \eps_0}{\sqrt{2}} \int_{\omega_{\eps_0,n}}|\psi_n(x)|^2dx >0.
\end{multline}
We readily obtain from \eqref{rr11} that
\begin{equation}\label{rr11k}
\liminf_{n \to +\infty} \frac{|\omega \cap [-\sqrt{2n+1},\sqrt{2n+1}]|}{\sqrt{n}}>0.
\end{equation}
Since the floor function satisfies 
\begin{equation*}
\forall R>0, \quad \lfloor R \rfloor \leq R < \lfloor R \rfloor+1,
\end{equation*}
it follows that for all $R \geq 2$,
\begin{equation}\label{floor_ineq1}
 \sqrt{2 n+1} \leq R < \sqrt{2 n+3},
\end{equation}
with $n=n(R)=\lfloor \frac{R^2-1}{2} \rfloor \in \nn$. We obtain that for all $R \geq 2$,
\begin{equation}\label{rr12k}
\frac{|\omega \cap [-R,R]|}{|[-R,R]|} \geq \frac{|\omega \cap [-\sqrt{2n+1},\sqrt{2n+1}]|}{2R} 
= \frac{\sqrt{n}}{2R} \frac{|\omega \cap [-\sqrt{2n+1},\sqrt{2n+1}]|}{\sqrt{n}}.
\end{equation}
Notice that \eqref{floor_ineq1} implies that
\begin{equation}\label{lim0}
\lim_{R \to +\infty} \frac{\sqrt{n}}{2R} = \frac{1}{2\sqrt{2}}.
\end{equation}
It follows from \eqref{rr11k}, \eqref{rr12k} and \eqref{lim0} that
\begin{equation*}
\liminf_{R \to +\infty} \frac{|\omega \cap [-R,R]|}{|[-R,R]|} \geq \frac{1}{2 \sqrt{2}} \liminf_{n \to +\infty} \frac{|\omega \cap [-\sqrt{2n+1},\sqrt{2n+1}]|}{\sqrt{n}} >0 .
\end{equation*}
It proves that assertion $(ii)$ implies assertion $(i)$.

Conversely, let $\omega$ be a measurable subset of $\rr$ verifying 
\begin{equation*}\label{rr12}
\liminf_{R \to +\infty} \frac{|\omega \cap [-R,R]|}{|[-R,R]|}>0
\end{equation*}
and define
\begin{equation*}
\tilde{\delta}= \liminf_{n \to +\infty} \frac{|\omega \cap [-\sqrt{2n+1}, \sqrt{2n+1}]|}{\sqrt{n}} >0.
\end{equation*}
In order to prove $(iii)$, it follows from \eqref{rr9b} that it is therefore sufficient to show that 
\begin{equation*}\label{rr13k}
\exists 0 <\eps <\frac{\pi}{2}, \quad \liminf_{n \to +\infty}\|\psi_n\|_{L^2(\omega_{\eps,n})}>0.
\end{equation*}
According to (\ref{rr9}) and (\ref{rr10}), it is then sufficient to show that 
\begin{equation}\label{rr13}
\exists 0 <\eps <\frac{\pi}{2}, \quad \liminf_{n \to +\infty}\|F_n\|_{L^2(\omega_{\eps,n})}>0.
\end{equation}
While using the substitution rule with $x=\sqrt{2n+1}\cos \theta$, we deduce from (\ref{e7}) and \eqref{rr9b} that for all $0<\eps <\frac{\pi}{2}$ and $n \geq 1$,
\begin{multline}\label{e12bb}
\int_{\omega_{\eps,n}}|F_n(x)|^2dx = \frac{1}{\pi}\sqrt{\frac{4n+2}{n}}\int_{\tilde{\omega}_{\eps,n}}\sin^2\Big(\Big(\frac{n}{2}+\frac{1}{4}\Big)\big(\sin(2\theta)-2\theta\big)+\frac{3\pi}{4}\Big)d\theta,
\end{multline}
with
\begin{equation}\label{e25}
\tilde{\omega}_{\eps,n}=\textrm{arccos}\Big(\frac{\omega_{\eps,n}}{\sqrt{2n+1}}\Big) \subset (\eps,\pi-\eps).
\end{equation}
According to (\ref{rr13}) and (\ref{e12bb}), 
it is then sufficient to show that 
\begin{equation}\label{rr13n}
\exists 0 <\eps <\frac{\pi}{2}, \quad \liminf_{n \to +\infty}\int_{\tilde{\omega}_{\eps,n}}\sin^2\Big(\Big(\frac{n}{2}+\frac{1}{4}\Big)\big(\sin(2\theta)-2\theta\big)+\frac{3\pi}{4}\Big)d\theta>0.
\end{equation}
Let $0<\eps_0 <\frac{\pi}{2}$ verifying 
\begin{equation*}\label{rrr1}
0<2\sqrt{2}(1-\cos \eps_0)<\tilde{\delta}.
\end{equation*}
By noticing from (\ref{rr9b}) that for all $0<\eps <\frac{\pi}{2}$ and $n \geq 1$,
\begin{equation*}\label{e13b}
\frac{|\omega_{\eps,n}|}{\sqrt{n}} \geq \frac{|\omega \cap [-\sqrt{2n+1},\sqrt{2n+1}]|}{\sqrt{n}}-\frac{2\sqrt{2n+1}}{\sqrt{n}}(1-\cos \eps),
\end{equation*}
we deduce that
\begin{equation*}\label{uu1}
\liminf_{n \to +\infty}\frac{|\omega_{\eps_0,n}|}{\sqrt{n}} \geq \tilde{\delta}-2\sqrt{2}(1-\cos \eps_0) >0.
\end{equation*}
It follows that there exist some positive constants $c_0>0$ and $n_0 \geq 1$ such that 
\begin{equation}\label{uu2}
\forall n \geq n_0, \quad |\omega_{\eps_0,n}| \geq c_0 \sqrt{n}.
\end{equation}
On the other hand, we notice from (\ref{e25}) and (\ref{uu2}) that for all $n \geq n_0$, 
\begin{equation}\label{uu3}
|\tilde{\omega}_{\eps_0,n}|=\int_{\omega_{\eps_0,n}}|\kappa_n'(x)|dx=\int_{\omega_{\eps_0,n}}\frac{dx}{\sqrt{(2n+1)-x^2}} \geq \frac{|\omega_{\eps_0,n}|}{\sqrt{2n+1}} \geq c_0 \sqrt{\frac{n}{2n+1}},
\end{equation}
with $\kappa_n(x)=\textrm{arccos}\big(\frac{x}{\sqrt{2n+1}}\big)$.
It follows from (\ref{uu3}) that there exists a positive constant $c_1>0$ such that 
\begin{equation}\label{uu4}
\forall n \geq n_0, \quad |\tilde{\omega}_{\eps_0,n}|\geq c_1 >0.
\end{equation}
By using anew the substitution rule with $t=2\theta-\sin(2\theta)$, we observe that 
\begin{multline}\label{uu5}
\int_{\tilde{\omega}_{\eps_0,n}}\sin^2\Big(\Big(\frac{n}{2}+\frac{1}{4}\Big)\big(\sin(2\theta)-2\theta\big)+\frac{3\pi}{4}\Big)d\theta\\
=\int_{\varphi(\tilde{\omega}_{\eps_0,n})}\frac{1}{\varphi'(\varphi^{-1}(t))}\sin^2\Big(\Big(\frac{n}{2}+\frac{1}{4}\Big)t-\frac{3\pi}{4}\Big)dt \geq c_2\int_{\varphi(\tilde{\omega}_{\eps_0,n})}\sin^2\Big(\Big(\frac{n}{2}+\frac{1}{4}\Big)t-\frac{3\pi}{4}\Big)dt,
\end{multline}
with 
$$0 < c_2=\inf_{t \in [2\eps_0-\sin(2\eps_0),2\pi-2\eps_0+\sin(2\eps_0)]}\frac{1}{\varphi'(\varphi^{-1}(t))}<+\infty,$$
where $\varphi(\theta)=2\theta-\sin(2\theta)$ is the increasing $C^{\infty}$-diffeomorphism from $[\eps_0,\pi-\eps_0]$ to $[2\eps_0-\sin(2\eps_0),2\pi-2\eps_0+\sin(2\eps_0)]$ already used above.
According to (\ref{rr13n}) and (\ref{uu5}), it therefore sufficient to prove that 
\begin{equation}\label{uu6}
\liminf_{n \to +\infty}\int_{\varphi(\tilde{\omega}_{\eps_0,n})}\sin^2\Big(\Big(\frac{n}{2}+\frac{1}{4}\Big)t-\frac{3\pi}{4}\Big)dt>0,
\end{equation}
where according to (\ref{uu4}), the subset $\varphi(\tilde{\omega}_{\eps_0,n}) \subset (2\eps_0-\sin(2\eps_0),2\pi-2\eps_0+\sin(2\eps_0))$ satisfies for all $n \geq n_0$,
\begin{equation}\label{uu7}
|\varphi(\tilde{\omega}_{\eps_0,n})|=\int_{\tilde{\omega}_{\eps_0,n}}|\varphi'(x)|dx \geq \Big(\inf_{[\eps_0,\pi-\eps_0]}|\varphi'|\Big) |\tilde{\omega}_{\eps_0,n}|\geq c_3,
\end{equation}
with 
$$0< c_3=c_1\inf_{[\eps_0,\pi-\eps_0]}|\varphi'|<+\infty.$$
By using anew the substitution rule with $x=(\frac{n}{2}+\frac{1}{4})t-\frac{3\pi}{4}$, we observe that 
\begin{equation}\label{uu8}
\int_{\varphi(\tilde{\omega}_{\eps_0,n})}\sin^2\Big(\Big(\frac{n}{2}+\frac{1}{4}\Big)t-\frac{3\pi}{4}\Big)dt=\frac{4}{2n+1}\int_{\Omega_{\eps_0,n}} \sin^2 x\  dx,
\end{equation}
where according to (\ref{e25}) and (\ref{uu7}), the subset 
\begin{equation*}\label{uu9}
\Omega_{\eps_0,n}=\Big(\frac{n}{2}+\frac{1}{4}\Big)\varphi(\tilde{\omega}_{\eps_0,n})-\frac{3\pi}{4} \subset \Big(\Big(\frac{n}{2}+\frac{1}{4}\Big)\varphi(\eps_0)-\frac{3\pi}{4},\Big(\frac{n}{2}+\frac{1}{4}\Big)\varphi(\pi-\eps_0)-\frac{3\pi}{4}\Big),
\end{equation*}
satisfies 
\begin{equation}\label{uu10}
\exists c_4>0, \forall n \geq n_0, \quad |\Omega_{\eps_0,n}| \geq c_4 n.
\end{equation}
According to (\ref{uu6}) and (\ref{uu8}), it is therefore sufficient to check that 
\begin{equation}\label{uu11}
\liminf_{n \to +\infty}\frac{1}{n}\int_{\Omega_{\eps_0, n}} \sin^2 x \ dx>0.
\end{equation}
Let $k \geq 1$ be an integer depending only on $\eps_0$ such that 
\begin{equation} \label{uu12}
 \forall n \geq n_0, \quad  \Omega_{\eps_0, n} \subset \Big(-kn\pi+\frac{\pi}{2}, kn \pi + \frac{\pi}{2}\Big).
\end{equation}
We denote $N_n = \#I_n$ the cardinality of the following set 
\begin{equation}\label{uu13}
I_n= \Big\{i \in [-kn, kn-1]\cap \mathbb{Z}:\ \Big|\Omega_{\eps_0, n} \cap \Big[i \pi+ \frac{\pi}{2}, \pi(i+1)+ \frac{\pi}{2}\Big]\Big| \geq \frac{c_4}{2k(\pi+1)} \Big\},
\end{equation}
with $n \geq n_0$.
It follows from \eqref{uu12} and \eqref{uu13} that for all $n \geq n_0$, 
\begin{align}\label{uu14}
|\Omega_{\eps_0,n}| &= \sum \limits_{i=-kn}^{kn-1}\Big|\Omega_{\eps_0, n} \cap \Big[i\pi+\frac{\pi}{2}, (i+1)\pi +\frac{\pi}{2}\Big] \Big| \\  \nonumber 
&=\sum \limits_{i \in I_n} \Big|\Omega_{\eps_0,n} \cap \Big[i\pi+\frac{\pi}{2}, (i+1)\pi +\frac{\pi}{2}\Big] \Big|+ \sum \limits_{\substack{i \notin I_n, \\ i \in [-kn,kn-1] \cap \mathbb{Z}}} \Big|\Omega_{\eps_0, n} \cap \Big[i\pi+\frac{\pi}{2}, (i+1)\pi +\frac{\pi}{2}\Big] \Big| \\  \nonumber 
&\leq \pi N_n + 2kn \frac{c_4}{2k(\pi+1)} = \pi N_n +\frac{nc_4}{\pi+1}.
\end{align}
We deduce from \eqref{uu10} and \eqref{uu14} that 
\begin{equation}\label{uu15}
\forall n \geq n_0, \quad N_n \geq \frac{n c_4}{\pi+1}.
\end{equation}
Let us now prove that the estimate \eqref{uu11} holds. By using the $\pi$-periodicity of the function $x \mapsto \sin^2 x$ and \eqref{uu12}, we obtain that for all $n \geq n_0$,
\begin{align}\label{uu16}
\frac{1}{n} \int_{\Omega_{\eps_0, n}} \sin^2x \ dx &= \frac{1}{n} \sum \limits_{i=-kn}^{kn-1} \int_{\Omega_{\eps_0, n} \cap [i\pi+\frac{\pi}{2}, (i+1) \pi+ \frac{\pi}{2}]} \sin^2 x \ dx
\\ \nonumber &\geq \frac{1}{n} \sum \limits_{i \in I_n} \int_{(\Omega_{\eps_0, n} -(i+1)\pi) \cap [-\frac{\pi}{2}, \frac{\pi}{2}]} \sin^2 x \ dx.
\end{align}
An application of Lemma~\ref{lemme_moyenne} in appendix together with \eqref{uu13}, \eqref{uu15} and \eqref{uu16} provide that for all $n \geq n_0$,
\begin{multline}\label{uu17}
\frac{1}{n} \int_{\Omega_{\eps_0, n}} \sin^2 x \ dx \geq \frac{1}{n} \sum \limits_{i \in I_n} \int_{-\frac{m_{i,n}}{2}}^{\frac{m_{i,n}}{2}} \sin^2 x \ dx \geq \frac{N_n}{n} \int_{-\frac{c_4}{4k(\pi+1)}}^{\frac{c_4}{4k(\pi+1)}} \sin^2 x \ dx \\
 \geq \frac{c_4}{\pi+1} \int_{-\frac{c_4}{4k(\pi+1)}}^{\frac{c_4}{4k(\pi+1)}} \sin^2 x \ dx >0,
\end{multline}
with
 $$m_{i,n} = \Big|(\Omega_{\eps_0, n} -(i+1)\pi) \cap \Big[-\frac{\pi}{2}, \frac{\pi}{2}\Big]\Big| = \Big|\Omega_{\eps_0, n} \cap \Big[i\pi+\frac{\pi}{2}, (i+1)\pi+\frac{\pi}{2}\Big]\Big|,$$
for all $n \geq n_0$ and $i \in [-kn,kn-1] \cap \mathbb{Z}$. This establishes \eqref{uu11} and concludes that assertion~$(i)$ implies assertion~$(iii)$. This ends the proof of  Proposition~\ref{prop}.

\section{Proof of Theorem~\ref{negative_result}}\label{proof_negative_result}
This section is devoted to the proof of Theorem~\ref{negative_result}. For any $\alpha=(\alpha_1,...,\alpha_d) \in (\nn \setminus \{0\})^d$ and $\eps \geq 0$, we consider the parallelepipeds
\begin{equation*}\label{pp1}
\mathcal{C}_{\alpha, \eps} = \prod_{j=1}^d [-\sqrt{2 \alpha_j+1}\cos \eps, \sqrt{2\alpha_j+1}\cos \eps].
\end{equation*}
Let $0<\eps <\frac{\pi}{2}$. We deduce from the one-dimensional Plancherel-Rotach formula \eqref{e4} that for all $\alpha \in (\nn \setminus \{0 \})^d$, $0<\eps < \frac{\pi}{2}$ and $x \in \mathcal{C}_{\alpha, \eps}$, 
\begin{equation}\label{decomp}
\Psi_{\alpha}(x)=F_{\alpha}(x)+ G_{\alpha}(x),
\end{equation}
with
\begin{equation}\label{decomp2}
F_{\alpha}(x)= \frac{2^{\frac{d}{4}}}{\pi^{\frac{d}{2}}\alpha^{\frac{1}{4}}} \prod_{i =1}^d \frac{\sin\big( \big(\frac{\alpha_i}{2}+\frac{1}{4} \big)(\sin(2 \theta_i)-2\theta_i)+\frac{3\pi}{4}\big)}{\sqrt{\sin \theta_i}}
\end{equation}
and
\begin{equation}\label{decomp3}
G_{\alpha}(x)= \frac{2^{\frac{d}{4}}}{\pi^{\frac{d}{2}}\alpha^{\frac{1}{4}}\prod_{i =1}^d \sqrt{\sin \theta_i}}R_{\alpha}(x),
\end{equation}
with $x=(\sqrt{2\alpha_1+1} \cos \theta_1,...,\sqrt{2\alpha_d+1} \cos \theta_d) \in \mathcal{C}_{\alpha, \eps}$.
The remainder terms $R_{\alpha}$ satisfy
\begin{equation}\label{esti1}
\exists C_{\eps}>0, \forall \alpha \in (\nn \setminus \{0\})^d, \forall x \in \mathcal{C}_{\alpha, \eps}, \quad |R_{\alpha}(x)| \leq C_{\eps} \sum_{i=1}^d \frac{1}{\alpha_i+1}.
\end{equation}
By using the preliminary step in the proof of Proposition~\ref{prop}, we obtain from (\ref{rr2}) that for all $0<\eps<\frac{\pi}{2}$,
\begin{align}\label{lim}
& \ \lim_{\alpha_{1} \to +\infty} ... \lim_{\alpha_{d} \to +\infty} \|\Psi_{\alpha} \|^2_{L^2(\mathcal{C}_{\alpha, \eps})}\\ \notag
 = & \ \lim_{\alpha_{1} \to +\infty} ... \lim_{\alpha_{d} \to +\infty} \prod_{j=1}^d \|\psi_{\alpha_j} \|^2_{L^2([-\sqrt{2\alpha_j+1}\cos \eps,\sqrt{2\alpha_j+1}\cos \eps])} \\ \notag
 = & \ \prod_{j=1}^d \lim_{\alpha_{j} \to +\infty} \|\psi_{\alpha_j} \|^2_{L^2([-\sqrt{2\alpha_j+1}\cos \eps,\sqrt{2\alpha_j+1}\cos \eps])} = \Big( 1- \frac{2\eps}{\pi} \Big)^{\frac{d}{2}}.
\end{align}
Let us assume that the fractional harmonic Schr\"odinger equation \eqref{harmonic_control} with $s \geq 1$ is exactly controllable from a measurable subset $\omega \subset \rr^d$ for some time $T>0$. The spectral gap condition in Proposition~\ref{hautus_test} is satisfied by the fractional harmonic oscillator $(-\Delta_x + |x|^2)^s$ since its spectrum is given by 
$$\sigma\big((-\Delta_x +|x|^2)^s\big) = \big\{(2n+d)^s: \  n \in \nn\big\}.$$
We can deduce from the infinite dimensional Hautus test \cite[Corollary~2.18]{miller} recalled in Proposition~\ref{hautus_test} that there exists $\delta>0$ such that 
\begin{equation} \label{sp}
\forall \alpha \in \nn^d, \quad \delta= \delta \|\Psi_{\alpha} \|_{L^2(\rr^d)} \leq \|\Psi_{\alpha} \|_{L^2(\omega)}.
\end{equation}
Let $0< \eps_0 < \frac{\pi}{2}$ such that
$$\delta -\Big(1-\Big( 1- \frac{2\eps_0}{\pi} \Big)^{\frac{d}{2}} \Big) >0.$$ 
By writing for all $\alpha \in \nn^d$,
\begin{equation*}\label{lim2}
\|\Psi_{\alpha} \|_{L^2(\omega)}^2 = \|\Psi_{\alpha} \|_{L^2(\omega\cap \mathcal{C}_{\alpha, \eps_0})}^2 + \|\Psi_{\alpha} \|_{L^2(\omega \setminus \mathcal{C}_{\alpha, \eps_0})}^2 \leq \|\Psi_{\alpha} \|_{L^2(\omega\cap \mathcal{C}_{\alpha, \eps_0})}^2+ 1-\|\Psi_{\alpha} \|_{L^2(\mathcal{C}_{\alpha, \eps_0})}^2,
\end{equation*}
we deduce from \eqref{lim} and \eqref{sp} that 
\begin{equation*}\label{lim3}
\liminf_{\alpha_{1} \to +\infty} ... \liminf_{\alpha_{d} \to +\infty} \|\Psi_{\alpha} \|_{L^2(\omega \cap \mathcal{C}_{\alpha, \eps_0})} \geq \sqrt{ \delta^2 -\Big(1-\Big( 1- \frac{2\eps_0}{\pi} \Big)^{\frac{d}{2}} \Big)} >0.
\end{equation*}
It follows from \eqref{decomp3} and \eqref{esti1} that 
\begin{equation*}
\forall \alpha \in (\nn \setminus \{0\})^d, \quad \|G_{\alpha} \|_{L^2(\mathcal{C}_{\alpha, \eps_0})}^2 \leq C_{\eps_0}^2 \frac{2^{\frac{d}{2}}}{ \alpha^{\frac{1}{2}}}\Big(\sum_{i=1}^d \frac{1}{\alpha_i+1} \Big)^2 \Big(\prod_{i=1}^d \sqrt{2 \alpha_i +1}\Big) \Big(1-\frac{2\eps_0}{\pi} \Big)^d
\end{equation*}
and
\begin{equation}\label{esti2}
\limsup_{\alpha_1 \to +\infty}... \limsup_{\alpha_d \to +\infty} \|G_{\alpha} \|_{L^2(\mathcal{C}_{\alpha, \eps_0})}=0.
\end{equation}
We can therefore deduce from \eqref{decomp} and \eqref{esti2} that
\begin{align}\label{lim4}
& \ \liminf_{\alpha_{1} \to +\infty} ... \liminf_{\alpha_{d} \to +\infty} \| F_{\alpha}\|_{L^2(\omega \cap \mathcal{C}_{\alpha, \eps_0})} \\ \notag
\geq & \ \liminf_{\alpha_{1} \to +\infty} ... \liminf_{\alpha_{d} \to +\infty}\big[\|\Psi_{\alpha} \|_{L^2(\omega \cap \mathcal{C}_{\alpha, \eps_0})} -\|G_{\alpha} \|_{L^2(\omega \cap\mathcal{C}_{\alpha, \eps_0})}\big] \\ \nonumber
 \geq & \  \liminf_{\alpha_{1} \to +\infty} ... \liminf_{\alpha_{d} \to +\infty}\big[\|\Psi_{\alpha} \|_{L^2(\omega \cap \mathcal{C}_{\alpha, \eps_0})} -\|G_{\alpha} \|_{L^2(\mathcal{C}_{\alpha, \eps_0})}\big] \\ \nonumber
 \geq & \ \liminf_{\alpha_{1} \to +\infty} ... \liminf_{\alpha_{d} \to +\infty}\|\Psi_{\alpha} \|_{L^2(\omega \cap \mathcal{C}_{\alpha, \eps_0})} >0.
\end{align}
On the other hand, we notice from \eqref{decomp2} that for all $\alpha \in (\nn \setminus \{0\})^d$,
\begin{equation*}
\forall x \in \mathcal{C}_{\alpha, \eps_0}, \quad |F_{\alpha}(x)|^2 \leq \frac{2^{\frac{d}{2}}}{\pi^d \alpha^{\frac{1}{2}}} \prod_{i=1}^d \frac{\sqrt{2\alpha_i+1}}{\sqrt{2\alpha_i+1-x_i^2}}
\end{equation*}
and this implies that for all $\alpha \in (\nn \setminus \{0\})^d$,
\begin{multline}\label{esti3}
\|F_{\alpha}\|^2_{L^2(\omega \cap \mathcal{C}_{\alpha, \eps_0})} \leq \frac{2^{\frac{d}{2}}}{\pi^d \alpha^{\frac{1}{2}}} \int_{\omega \cap \mathcal{C}_{\alpha, \eps_0}} \prod_{i =1}^d \frac{\sqrt{2\alpha_i+1}}{\sqrt{2\alpha_i+1-x_i^2}} dx \\ 
\leq \frac{2^{\frac{d}{2}}}{\pi^d\sin^d \eps_0} \frac{|\omega \cap [-\sqrt{2 \alpha_1+1}, \sqrt{2\alpha_1+1}]\times...\times [-\sqrt{2 \alpha_d+1}, \sqrt{2\alpha_d+1}]|}{\alpha^{\frac{1}{2}}}.
\end{multline}
We deduce from \eqref{lim4} and \eqref{esti3} that
\begin{equation}\label{m1}
\liminf_{\alpha_{1} \to +\infty} ... \liminf_{\alpha_{d} \to +\infty} \frac{|\omega \cap [-\sqrt{2 \alpha_1+1}, \sqrt{2\alpha_1+1}]\times ... \times [-\sqrt{2 \alpha_d+1}, \sqrt{2\alpha_d+1}]|}{|[-\sqrt{2 \alpha_1+1}, \sqrt{2\alpha_1+1}]\times ... \times [-\sqrt{2 \alpha_d+1}, \sqrt{2\alpha_d+1}]|} >0.
\end{equation}
It follows from (\ref{floor_ineq1}) that for all $(R_1,...,R_d) \in [2, +\infty)^d$,
\begin{equation}\label{floor_ineq2}
\forall 1 \leq i \leq d, \quad \sqrt{2 n_i+1} \leq R_i < \sqrt{2 n_i+3},
\end{equation}
where $n_i=n_i(R_i)= \lfloor \frac{R_i^2-1}{2} \rfloor \in \nn$ for all $1 \leq i \leq d$. We obtain that
\begin{align}\label{m2}
& \  \frac{|\omega \cap [-R_1, R_1]\times ... \times [-R_d, R_d]|}{|[-R_1, R_1]\times ... \times [-R_d, R_d]|} \\  \notag
\geq &\  \frac{|\omega \cap [-\sqrt{2 n_1+1}, \sqrt{2n_1 +1}]\times ... \times [-\sqrt{2 n_d+1}, \sqrt{2 n_d+1}]|}{|[-R_1, R_1]\times ... \times [-R_d, R_d]|} \\ \notag
= & \ \Big(\prod_{i=1}^d \frac{\sqrt{2n_i+1}}{R_i}\Big) \frac{|\omega \cap [-\sqrt{2 n_1+1}, \sqrt{2 n_1+1}]\times ... \times [-\sqrt{2 n_d+1}, \sqrt{2 n_d+1}]|}{|[-\sqrt{2 n_1+1}, \sqrt{2 n_1+1}]\times ... \times [-\sqrt{2 n_d+1}, \sqrt{2 n_d+1}]|} .
\end{align}
We observe that \eqref{floor_ineq2} implies that
\begin{equation}\label{lim5}
\forall 1 \leq i \leq d, \quad \lim_{R_i \to +\infty} \frac{\sqrt{2n_i+1}}{R_i} = 1.
\end{equation}
It follows from \eqref{m1}, \eqref{m2} and \eqref{lim5} that 
\begin{multline*}
\liminf_{R_{1} \to +\infty} ... \liminf_{R_{d} \to +\infty} \frac{|\omega \cap [-R_1, R_1]\times ... \times [-R_d, R_d]|}{|[-R_1, R_1]\times ... \times [-R_d, R_d]|} \\ \geq 
\liminf_{\alpha_{1} \to +\infty} ... \liminf_{\alpha_{d} \to +\infty} \frac{|\omega \cap [-\sqrt{2 \alpha_1+1}, \sqrt{2\alpha_1+1}]\times ... \times [-\sqrt{2 \alpha_d+1}, \sqrt{2\alpha_d+1}]|}{|[-\sqrt{2 \alpha_1+1}, \sqrt{2\alpha_1+1}]\times ... \times [-\sqrt{2 \alpha_d+1}, \sqrt{2\alpha_d+1}]|} >0.
\end{multline*}
It establishes \eqref{necessary_cond} with $A= I_d$.
Let $A \in O(\rr^d)$. By using the invariance of the Weyl quantization under affine symplectic transformations (see e.g.~\cite[Theorem~18.5.9]{hormander}), we notice that
\begin{equation}\label{commute}
R_A \circ \mathcal{H} = \mathcal{H} \circ R_A,
\end{equation}
with $\mathcal{H}=-\Delta_x+|x|^2$, where $R_A$ is the unitary operator on $L^2(\rr^d)$ defined by $R_A f = f(A \, \cdot )$ for $f \in L^2(\rr^d)$.
The identity (\ref{commute}) implies in particular that $(R_A\Psi_{\alpha})_{\alpha \in \nn^d}$ is a Hilbert basis of $L^2(\rr^d)$ satisfying
$$\forall \alpha \in \nn^d, \quad \mathcal{H}(R_A\Psi_{\alpha})=(2|\alpha|+d)R_A\Psi_{\alpha},$$
when $A \in O(\rr^d)$. By using that the operator 
$$R_A : \textrm{Span}\big((\Psi_{\alpha})_{\alpha \in \nn^d, |\alpha|=N}\big) \rightarrow \textrm{Span}\big((\Psi_{\alpha})_{\alpha \in \nn^d, |\alpha|=N}\big),$$
is a unitary operator for any $N \in \nn$, we notice that for all $s \geq 1$ and $f \in \mathcal{D}(\mathcal{H}^s)$,
\begin{align*}
(R_A \circ \mathcal{H}^s)f=& \ (R_A \circ \mathcal{H}^s)\Big(\sum_{\alpha \in \nn^d}\langle f,R_{A^{-1}}\Psi_{\alpha}\rangle_{L^2(\rr^d)}R_{A^{-1}}\Psi_{\alpha}\Big)\\
= & \ R_A\Big(\sum_{\alpha \in \nn^d}(2|\alpha|+d)^s\langle f,R_{A^{-1}}\Psi_{\alpha}\rangle_{L^2(\rr^d)}R_{A^{-1}}\Psi_{\alpha}\Big)\\
=& \ \sum_{\alpha \in \nn^d}(2|\alpha|+d)^s\langle R_Af,\Psi_{\alpha}\rangle_{L^2(\rr^d)}\Psi_{\alpha}=(\mathcal{H}^s \circ R_A)f,
\end{align*}
that is
\begin{equation}\label{gj3}
\forall A \in O(\rr^d), \forall s \geq 1, \quad R_A \circ \mathcal{H}^s=\mathcal{H}^s \circ R_A.
\end{equation}
Let $f_0, f_T \in L^2(\rr^d)$. By assumption, there exists a control function $u \in L^2([0,T] \times \rr^d)$ supported in $\rr_+ \times \omega$ such that the semigroup solution $g : \rr_+ \times \rr^d \rightarrow \rr$ of the system \eqref{harmonic_control} associated to the initial condition $f_0(A \, \cdot)$ satisfies $g(T, \cdot)=f_T(A \, \cdot)$. We deduce from \eqref{gj3} that the function $f$ defined on $\rr_+ \times \rr^d$ by 
\begin{equation*}
\forall t \in \rr_+, \forall x \in \rr^d, \quad f(t,x)= g(t,A^{-1}x),
\end{equation*}
is the semigroup solution of \eqref{harmonic_control} associated to the initial condition $f_0$ and the control function 
$v(t,x)=u(t,A^{-1}x) \in L^2([0,T] \times \rr^d)$ supported in $\rr_+ \times A(\omega)$; and satisfies $f(T,\cdot)=f_T(\cdot)$. It implies that the harmonic Schr\"odinger equation \eqref{harmonic_control} is exactly controllable from the control set $A(\omega)$. By using the first part of this proof with $\omega$ replaced by $A(\omega)$, we obtain that 
\begin{equation*}
\liminf_{R_{1} \to +\infty} ... \liminf_{R_{d} \to +\infty} \frac{|A(\omega) \cap [-R_1, R_1]\times ... \times [-R_d, R_d]|}{|[-R_1, R_1]\times ... \times [-R_d, R_d]|} >0.
\end{equation*}
This ends the proof of Theorem~\ref{negative_result}.

\section{Appendix}
This appendix contains the proof of two instrumental results. The first part is devoted to the proof of a technical lemma used in the proof of Proposition~\ref{prop}. The purpose of the second part is to recall the proof of Miller's result stated in Proposition~\ref{ineqspec}.
\subsection{A technical lemma}
The following technical result is used in the proof of Proposition~\ref{prop}:

\medskip

\begin{lemma}\label{lemme_moyenne}
If $A$ is a measurable subset of $[-\frac{\pi}{2}, \frac{\pi}{2}]$ and $\delta=|A|$ denotes its Lebesgue measure then
\begin{equation*}
\int_A \sin^2 x \ dx \geq \int_{-\frac{\delta}{2}}^{\frac{\delta}{2}} \sin^2 x \ dx.
\end{equation*}
\end{lemma}

\medskip

\begin{proof}
By using that the function $x \mapsto \sin^2 x$ is even and increasing on $[0,\frac{\pi}{2}]$, we observe that
\begin{align}\label{lm1}
& \ \int_A \sin^2 x \ dx = \int_{A \cap [-\frac{\delta}{2}, \frac{\delta}{2} ]} \sin^2 x \ dx + \int_{A \cap ([-\frac{\pi}{2},\frac{\pi}{2}] \setminus [-\frac{\delta}{2}, \frac{\delta}{2}])} \sin^2 x \ dx \\ \notag
\geq & \ \int_{A \cap [-\frac{\delta}{2}, \frac{\delta}{2} ]} \sin^2 x \ dx+ \Big|A \cap \Big(\Big[-\frac{\pi}{2},\frac{\pi}{2}\Big] \setminus \Big[-\frac{\delta}{2}, \frac{\delta}{2}\Big]\Big) \Big| \sin^2\Big(\frac{\delta}{2} \Big) \\ \notag
 \geq & \ \int_{A \cap [-\frac{\delta}{2}, \frac{\delta}{2} ]} \sin^2 x \ dx+ \frac{\big|A \cap \big([-\frac{\pi}{2},\frac{\pi}{2}] \setminus [-\frac{\delta}{2}, \frac{\delta}{2}]\big) \big|}{\big| [-\frac{\delta}{2}, \frac{\delta}{2}] \cap \big([-\frac{\pi}{2},\frac{\pi}{2}]\setminus A\big)  \big|}\int_{[-\frac{\delta}{2}, \frac{\delta}{2}] \cap ([-\frac{\pi}{2},\frac{\pi}{2}]\setminus A)} \sin^2 x \ dx.
\end{align}
By using that $\delta=|A|=|[-\frac{\delta}{2}, \frac{\delta}{2}]|$, we notice that
\begin{multline}\label{lm2}
\Big|A \cap \Big(\Big[-\frac{\pi}{2},\frac{\pi}{2}\Big] \setminus \Big[-\frac{\delta}{2}, \frac{\delta}{2}\Big]\Big) \Big| = |A| -\Big|A \cap \Big[-\frac{\delta}{2}, \frac{\delta}{2}\Big]\Big|  \\ 
=  \Big|\Big[-\frac{\delta}{2}, \frac{\delta}{2}\Big]\Big|-\Big|A \cap \Big[-\frac{\delta}{2}, \frac{\delta}{2}\Big]\Big| = \Big| \Big[-\frac{\delta}{2}, \frac{\delta}{2}\Big] \cap \Big(\Big[-\frac{\pi}{2},\frac{\pi}{2}\Big]\setminus A\Big) \Big|.
\end{multline}
It follows from \eqref{lm1} and \eqref{lm2} that
\begin{equation*}
\int_A \sin^2 x \ dx \geq \int_{A \cap [-\frac{\delta}{2}, \frac{\delta}{2} ]} \sin^2 x \ dx+ \int_{[-\frac{\delta}{2}, \frac{\delta}{2}] \cap ([-\frac{\pi}{2},\frac{\pi}{2}]\setminus A)} \sin^2 x \ dx = \int_{-\frac{\delta}{2}}^{\frac{\delta}{2}} \sin^2 x \ dx.
\end{equation*}
This ends the proof of Lemma~\ref{lemme_moyenne}.
\end{proof}

\subsection{Spectral inequalities and exact controllability}\label{proof_ineqspec}
This section is devoted to recall the proof of Miller's result \cite[Corollary~2.17]{miller} stated in Proposition~\ref{ineqspec} which provides necessary and sufficient spectral estimates for the observability of system \eqref{obs01} to hold. The proof of Proposition~\ref{ineqspec} is based on another Miller's result \cite[Theorem~2.4]{miller} characterizing the observability with resolvent conditions which is recalled below and whose proof is omitted:

\medskip

\begin{proposition}[Miller {\cite[Theorem~2.4]{miller}}]\label{prop_resolv_condition}
Let $(A,\mathcal{D}(A))$ be a selfadjoint operator on $L^2(\rr^d)$, which is the infinitesimal generator of a strongly continuous group $(e^{itA})_{t \in \rr}$ on $L^2(\rr^d)$. 
The system \eqref{obs01} is exactly observable from a measurable subset $\omega \subset \rr^d$ if and only if there exist some positive constants $M>0$ and $m>0$ such that 
\begin{equation}\label{resolvent_condition}
\forall f \in \mathcal{D}(A), \, \forall \lambda \in \rr, \quad \|f\|^2_{L^2(\rr^d)} \leq M \| (A- \lambda) f \|^2_{L^2(\rr^d)} + m \|f\|^2_{L^2(\omega)}.
\end{equation}
When condition \eqref{resolvent_condition} is satisfied, exact observability holds in any time $T > \pi \sqrt{M}$.
\end{proposition}

\medskip

We consider the system \eqref{obs01}. If this system is exactly observable from a measurable subset $\omega \subset \rr^d$ at some time $T>0$, Proposition~\ref{prop_resolv_condition} proves that there exist some positive constants $M>0$ and $m>0$ such that the resolvent estimate \eqref{resolvent_condition} holds.
Let $\lambda \in \rr$, $0< D < \frac{1}{M}$ and $f \in \un_{\{|A-\lambda| \leq \sqrt{D}\}}\big(\mathcal{D}(A)\big)$. The functional calculus shows that 
\begin{equation}\label{fc1}
\|(A-\lambda)f \|_{L^2(\rr^d)} \leq \sqrt{D} \|f \|_{L^2(\rr^d)}.
\end{equation}
It follows from \eqref{resolvent_condition}  and \eqref{fc1} that 
\begin{equation}\label{rc1}
\forall \lambda \in \rr, \forall f \in\un_{\{|A-\lambda| \leq \sqrt{D}\}}\big(\mathcal{D}(A)\big), \quad \|f\|^2_{L^2(\rr^d)} \leq M D \|f\|^2_{L^2(\rr^d)} +m \|f\|^2_{L^2(\omega)},
\end{equation}
that is
\begin{equation*}
\forall \lambda \in \rr, \forall f \in\un_{\{|A-\lambda| \leq \sqrt{D}\}}\big(\mathcal{D}(A)\big), \quad \|f\|^2_{L^2(\rr^d)} \leq \frac{m}{1-MD} \|f\|^2_{L^2(\omega)}.
\end{equation*}
It establishes the spectral estimates \eqref{Spectral_ineq} for all $f \in \un_{\{|A-\lambda| \leq \sqrt{D}\}}\big(L^2(\rr^d)\big)$ since the domain $\mathcal{D}(A)$ is dense in $L^2(\rr^d)$.

Conversely, let us assume that there exist some positive constants $D>0$ and $k>0$ such that the spectral estimates \eqref{Spectral_ineq} holds. Let $\lambda \in \rr$ and $f \in \mathcal{D}(A)$. With $f_{\lambda} = \un_{\{|A-\lambda| \leq \sqrt{D} \}} f$ and $f_{\perp}=f-f_{\lambda}$, we deduce from \eqref{Spectral_ineq} and the functional calculus that for all $\lambda \in \rr$ and $f \in \mathcal{D}(A)$,
\begin{multline}\label{rc2}
\|f\|^2_{L^2(\rr^d)} = \|f_{\lambda} \|^2_{L^2(\rr^d)} + \|f_{\perp}\|^2_{L^2(\rr^d)} \leq k \|f_{\lambda}\|^2_{L^2(\omega)} + \frac{1}{D} \|(A- \lambda) f_{\perp}\|^2_{L^2(\rr^d)}  \\ \leq k \|f_{\lambda}\|^2_{L^2(\omega)} + \frac{1}{D} \|(A- \lambda) f \|^2_{L^2(\rr^d)}.
\end{multline}
Let $T> \pi \sqrt{\frac{1+k}{D}}$ and $\eps >0$ such that 
$$T^2 >\pi^2 \frac{1+(1+\eps^2)k}{D}.$$ 
The functional calculus shows that for all $\lambda \in \rr$ and $f \in \mathcal{D}(A)$,
\begin{multline}\label{rc3}
 \|f_{\lambda}\|^2_{L^2(\omega)} =\|f-f_{\perp}\|^2_{L^2(\omega)} \leq (1+\eps^{-2}) \|f\|^2_{L^2(\omega)} + (1+\eps^2) \|f_{\perp} \|^2_{L^2(\omega)} \\
\leq (1+\eps^{-2}) \|f\|^2_{L^2(\omega)} + (1+\eps^2) \|f_{\perp} \|^2_{L^2(\rr^d)}   \leq (1+\eps^{-2}) \|f\|^2_{L^2(\omega)} + \frac{1+\eps^2}{D} \|(A-\lambda) f_{\perp} \|^2_{L^2(\rr^d)} \\
\leq  (1+\eps^{-2}) \|f\|^2_{L^2(\omega)} + \frac{1+\eps^2}{D} \|(A-\lambda) f \|^2_{L^2(\rr^d)}.
\end{multline}
 It follows from \eqref{rc2} and \eqref{rc3} that 
\begin{equation*}
\forall \lambda \in \rr, \forall f \in \mathcal{D}(A), \quad \|f\|^2_{L^2(\rr^d)}  \leq k(1+\eps^{-2}) \|f\|^2_{L^2(\omega)} + \frac{1+k(1+\eps^2)}{D} \|(A- \lambda) f \|^2_{L^2(\rr^d)}.
\end{equation*}
We can deduce from Proposition~\ref{prop_resolv_condition} that the system \eqref{obs01} is exactly observable from $\omega$ in time $T$.

\end{document}